\documentstyle[12pt]{article}  
\def\sq{\hbox {\rlap{$\sqcap$}$\sqcup$}}
\def\sq{\hbox {\rlap{$\sqcap$}$\sqcup$}}
\def\R{ {\rm R \kern -.31cm I \kern .15cm}}
\def\C{ {\rm C \kern -.15cm \vrule width.5pt \kern .12cm}}
\def\Z{ {\rm Z \kern -.27cm \angle \kern .02cm}}
\def\N{ {\rm N \kern -.26cm \vrule width.4pt \kern .10cm}}
\def\1{{\rm 1\mskip-4.5mu l} }
\def\lsim{\raise0.3ex\hbox{$<$\kern-0.75em\raise-1.1ex\hbox{$\sim$}}}
\def\gsim{\raise0.3ex\hbox{$>$\kern-0.75em\raise-1.1ex\hbox{$\sim$}}}
\def\noi{\noindent}

\def\beq{\begin{equation}}   \def\eeq{\end{equation}}
\def\bea{\begin{eqnarray}}  \def\eea{\end{eqnarray}}
\def\nn{\nonumber}
\def\noi{\noindent}
\def\beeq{\begin{eqnarray}} \def\eeeq{\end{eqnarray}}
\newcommand\mysection{\setcounter{equation}{0}\section}

\newcounter{hran}

\begin{document} 
\centerline{\Large\bf Long Range Scattering and Modified } 
 \vskip 3 truemm \centerline{\Large\bf  Wave Operators for the Wave-Schr\"odinger
System II\footnote{Work supported in part by NATO Collaborative Linkage Grant 979341}} 

\vskip 0.5 truecm

\centerline{\bf J. Ginibre}
\centerline{Laboratoire de Physique Th\'eorique\footnote{Unit\'e Mixte de
Recherche (CNRS) UMR 8627}}  \centerline{Universit\'e de Paris XI, B\^atiment
210, F-91405 ORSAY Cedex, France}
\vskip 3 truemm
\centerline{\bf G. Velo}
\centerline{Dipartimento di Fisica, Universit\`a di Bologna}  \centerline{and INFN, Sezione di
Bologna, Italy}

\vskip 1 truecm

\begin{abstract}
We continue the study of scattering theory for the system consisting of a Schr\"odinger equation and a
wave equation with a Yukawa type coupling in space dimension 3. In a previous paper, we proved the
existence of modified wave operators for that system with no size restriction on the data and we
determined the asymptotic behaviour in time of solutions in the range of the wave operators, under a
support condition on the asymptotic state required by the different propagation properties of the wave
and Schr\"odinger equations. Here we eliminate that condition by using an improved asymptotic form for
the solutions.
\end{abstract}

\vskip 3 truecm
\noi AMS Classification : Primary 35P25. Secondary 35B40, 35Q40, 81U99.  \par
\noi Key words : Long range scattering, modified wave operators, Wave-Schr\"odinger system. 
\vskip 1 truecm

\noindent LPT Orsay 03-07 \par
\noindent February 2003 \par

\newpage
\pagestyle{plain}
\baselineskip 18pt

\mysection{Introduction}
\hspace*{\parindent}
This paper is a sequel to a previous paper with the same title (\cite{1r}, hereafter referred to as I)
where we studied the theory of scattering and proved the existence of modified wave operators for the
Wave-Schr\"odinger (WS) system
$$\hskip 3 truecm \left \{ \begin{array}{l} i\partial_t u = - (1/2) \Delta u - Au \hskip 6.5 truecm (1.1)
\\ \\ \sq A = |u|^2 \hskip 8.9 truecm (1.2)
\end{array}\right .$$

\noi where $u$ and $A$ are respectively a complex valued and a real valued function defined in space
time ${I\hskip-1truemm R}^{3+1}$. We refer to the introduction of I for general background and
references and we give here only a general overview of the problem. \par

The main result of I was the construction of modified wave operators for the WS system, with no size
restriction on the solutions. That construction basically consists in solving the Cauchy problem for the
WS system with infinite initial time, namely in constructing solutions with prescribed asymptotic
behaviour at infinity in time. That asymptotic behaviour is imposed in the form of suitable approximate
solutions of the WS system. One then looks for exact solutions, the difference of which with the given
approximate ones tends to zero at infinity in time in a suitable sense, more precisely in suitable
norms. The approximate solutions are obtained as low order iterates in an iterative resolution scheme of
the WS system. In I we used second order iterates. They are parametrized by data $(u_+, A_+, \dot{A}_+)$
which play the role of (actually are in simpler cases) initial data at time zero. Those data constitute
the asymptotic state for the actual solution. \par

An inherent difficulty of the WS system is the difference of propagation properties of the wave
equation and of the Schr\"odinger equation. Because of that difficulty, we had to impose in I a
support condition on the Fourier transform $Fu_+$ of the Schr\"odinger asymptotic state $u_+$, saying
in effect that $Fu_+$ vanishes in a neighborhood of the unit sphere, so that $u_+$ generates a
solution of the free Schr\"odinger equation which is asymptotically small in a neighborhood of the
light cone. Such a support condition is unpleasant because it cannot be satisfied on a dense subspace
of any reasonable space where one hopes to solve the problem, typically with $u$ in $FH^k$ for $k > 1/2$
($H^k$ is the standard $L^2$ based Sobolev space). \par

The theory of scattering and the existence of modified wave operators can also be studied for various
equations and systems including the WS system by a method simpler than that of I, proposed earlier by
Ozawa \cite{5r}. Contrary to that of I, that method is restricted to the case of small data and small
solutions. It has been applied to various systems, in particular to the Klein-Gordon-Schr\"odinger
(KGS) system in dimension 2, which is fairly similar to the WS system in dimension 3 from the point of
view of scattering \cite{6r}. Similar propagation difficulties also appear for that system, thereby
again requiring a support condition on $Fu_+$ in the treatment given in \cite{6r}. \par

A progress on that problem was made recently by Shimomura \cite{7r} \cite{8r} who was able to remove
the previous support condition in the construction of the modified wave operators by the Ozawa method
in the case of the KGS system in dimension 2 \cite{7r} and of the WS system in dimension 3 \cite{8r}.
The key of that progress consists in using an improved asymptotic form for the Schr\"odinger
function, obtained by adding a term depending on $(A_+, \dot{A}_+)$ which partly cancels the
contribution of the asymptotic field for $A$ in the Schr\"odinger equation. \par

Although the method used in I is more complicated than the Ozawa method (so as to accomodate
arbitrarily large data and solutions), it turns out that the improved asymptotic form of $u$ used in
\cite{8r} can be transposed into the framework of the method of I, thereby allowing to remove the
support condition on $Fu_+$ assumed in I, at least in a restricted interval of values of the
parameters defining the regularity of the solutions. The purpose of the present paper is to implement
that improvement, namely to rederive the main results of I without assuming the support condition
used in I, by using the improved asymptotic form of the solution inspired by that of \cite{8r}. \par

In the remaining part of this introduction, we shall briefly review the method used in I in the modified
form used in the present paper. We refer to Section 2 of I for a more detailed exposition. The main result
of this paper will be stated in semi heuristic terms at the end of this introduction. The first step in
that method consists in eliminating the wave equation (1.2) by solving it for $A$ and substituting the
result into the Schr\"odinger equation, which then becomes both non linear and non local in time. One then
parametrizes the Schr\"odinger function $u$ in terms of an amplitude $w$ and a phase $\varphi$ and one
replaces the Schr\"odinger equation by an auxiliary system consisting of a transport equation for the
amplitude and a Hamilton-Jacobi equation for the phase. One solves the Cauchy problem at infinity, namely
with prescribed asymptotic behaviour, for the auxiliary system, and one finally reconstructs the solution
of the original WS system from that of the auxiliary system. We now proceed to the technical details. We
restrict our attention to positive time, actually to $t \geq 1$. \par

We first eliminate the wave equation. We define 
$$\omega = (- \Delta)^{1/2} \quad, \quad K(t) = \omega^{-1} \sin \omega t \quad , \quad \dot{K}(t) =
\cos \omega t$$    

\noi and we replace (1.2) by
$$A = A_0 + A_1 (|u|^2) \eqno(1.3)$$

\noi where
$$A_0 = \dot{K}(t) A_+ + K(t) \dot{A}_+ \ , \eqno(1.4)$$
$$A_1 (|u|^2) = - \int_t^{\infty} dt'\ K(t-t') |u(t')|^2 \ . \eqno(1.5)$$

\noi Here $A_0$ is a solution of the free wave equation with initial data $(A_+, \dot{A}_+)$ at time $t
= 0$. The pair $(A_+, \dot{A}_+)$ is the asymptotic state for $A$. \par

We next perform the change of variables mentioned above on $u$. The unitary
group 
$$U(t) = \exp (i(t/2)\Delta ) \eqno(1.6)$$

\noi which solves the free Schr\"odinger equation can be written as
$$U(t) = M(t) \ D(t) \ F \ M(t) \eqno(1.7)$$

\noi where $M(t)$ is the operator of multiplication by the function
$$M(t) = \exp  \left ( i x^2/2t \right ) \ , \eqno(1.8)$$

\noi $F$ is the Fourier transform and $D(t)$ is the dilation operator
$$D(t) = (it)^{-n/2} \ D_0(t) \eqno(1.9)$$

\noi where
$$\left ( D_0(t) f\right )(x) = f(x/t) \ . \eqno(1.10)$$

\noi We parametrize $u$ in terms of an amplitude $w$ and of a real phase $\varphi$ as 
$$u(t) = M(t) \ D(t) \exp [- i \varphi (t)] w(t) \ . \eqno(1.11)$$

\noi Substituting (1.11) into (1.1) yields an evolution equation for $(w , \varphi )$,
namely  
$$\left \{ i \partial_t + (2t^2)^{-1} \Delta - i(2t^2)^{-1} (2 \nabla \varphi \cdot \nabla +
\Delta \varphi ) + t^{-1} B + \partial_t \varphi - (2t^2)^{-1} |\nabla \varphi |^2 \right \} w
= 0 \eqno(1.12)$$

\noi where we have expressed $A$ in terms of a new function $B$ by
$$A = t^{-1} \ D_0\ B \ . \eqno(1.13)$$

\noi Corresponding to the decomposition (1.3) of $A$, we decompose
$$B = B_0 + B_1 (w , w) \eqno(1.14)$$

\noi where $A_0 = t^{-1} D_0B_0$ and $A_1 = t^{-1}D_0B_1$. One computes easily
$$B_1(w_1,w_2) = \int_1^{\infty} d \nu \ \nu^{-3} \ \omega^{-1} \sin ((\nu - 1) \omega)
D_0(\nu) ({\rm Re} \ \bar{w}_1w_2)(\nu t) \ .\eqno(1.15)$$

\indent At this point, we have only one evolution equation (1.12) for two
functions $(w, \varphi )$. We arbitrarily impose a second equation, namely a Hamilton-Jacobi (or
eikonal) equation for the phase $\varphi$, thereby splitting the equation (1.12) into a system of two
equations, the other one of which being a transport type equation for the amplitude $w$.
For that purpose, we split $B_0$ and $B_1$ into a long range and short range parts as follows. Let
$\chi \in {\cal C}^{\infty} ({I\hskip-1truemm R}^{3}, {I\hskip-1truemm R})$, $0 \leq \chi \leq 1$, $\chi
(\xi) = 1$ for $|\xi| \leq 1$, $\chi (\xi) = 0$ for $|\xi | \geq 2$ and let $0 < \beta_0$, $\beta < 1$.
We define
$$B_0 = B_{0L} + B_{0S} \qquad , \qquad B_1 = B_L + B_S \eqno(1.16)$$

\noi where
$$\left \{ \begin{array}{l}  FB_{0L}(t, \xi ) = \chi (\xi t^{-\beta_0})
FB_0(t, \xi ) \ , \\ \\
 FB_L(t, \xi ) = \chi (\xi  t^{-\beta})
FB_1(t, \xi ) \ .   \end{array}
\right . \eqno(1.17)$$ 

\noi The splitting (1.16) (1.17) differs from that made in I in two respects. First and more important
is the fact that we perform that splitting both on $B_0$ and on $B_1$, whereas in I it was done only on
$B_1$. Second, we use here a smooth cut-off $\chi$ instead of a sharp one. The smooth cut-off is
actually needed only for $B_0$. For $\beta = \beta_0$, the splitting is the same for $B_0$ and $B_1$
and can therefore be performed on $B$ without any reference to the asymptotic state $(A_+,
\dot{A}_+)$. The parameters $\beta_0$ and $\beta$ will have to satisfy various conditions which will
appear later, all of them compatible with $\beta = \beta_0 = 1/3$. \par

We split the equation (1.12) into the following system of two equations.
$$\left \{ \begin{array}{l} \partial_t w = i(2t^2)^{-1} \Delta w + t^{-2} Q (\nabla \varphi , w)
+ it^{-1} (B_{0S} + B_S(w, w)) w \\ \\
 \partial_t \varphi = (2t^2)^{-1} |\nabla \varphi |^2 - t^{-1} \ B_{0L}- t^{-1} \ B_L(w,w)
\end{array} \right . \eqno(1.18)$$ 
\noi where we have defined
$$Q(s, w) = s \cdot \nabla w + (1/2) (\nabla \cdot s) w
\eqno(1.19)$$

\noi for any vector field $s$. The first equation of (1.18) is the transport type equation for
the amplitude $w$, while the second one is the Hamilton-Jacobi type equation for the phase $\varphi$.
Since the right-hand sides of (1.18) contain $\varphi$ only through its gradient, we can obtain from
(1.18) a closed system for $w$ and $s = \nabla \varphi$ by taking the gradient of the second equation,
namely  
$$\left \{ \begin{array}{l} \partial_t w = i(2t^2)^{-1} \Delta w + t^{-2} Q (s , w)
+ it^{-1} (B_{0S} + B_S(w, w)) w \\ \\
 \partial_t s = t^{-2} s \cdot \nabla s - t^{-1} \nabla  B_{0L} - t^{-1} \nabla  B_L(w,w) \ .\end{array}
\right .  \eqno(1.20)$$ 

 Once the system (1.20) is solved for $(w, s)$, one recovers $\varphi$ easily by 
integrating the second equation of (1.18) over time. The system (1.20) will be referred to as the
auxiliary system.\par

The construction of the modified wave operators follows the same pattern as in I. The first task is to
construct solutions of the auxiliary system (1.20) with suitably prescribed asymptotic behaviour at
infinity, and in particular with $w(t)$ tending to a limit $w_+ = Fu_+$ as $t \to \infty$.  That
asymptotic behaviour is imposed in the form of a suitably chosen pair $(W, \phi )$ and therefore 
$(W, S)$ with $S = \nabla \phi$ with $W(t)$ tending to $w_+$ as $t \to \infty$. For fixed $(W,S)$, we
make a change of variables in the system (1.18) from $(w, \varphi )$ to $(q, \psi )$ defined by 
$$(q, \psi) = (w, \varphi ) - (W, \phi )  \eqno(1.21)$$
\noi or equivalently a change of variables in the system (1.20) from $(w, s)$ to $(q, \sigma )$ defined
by   
$$(q, \sigma) = (w, s ) - (W, S) \ , \eqno(1.22)$$

\noi and instead of looking for a solution $(w, s)$ of the system (1.20) with $(w, s)$ behaving
asymptotically as $(W, S)$, we look for a solution $(q, \sigma )$ of the transformed system with $(q,
\sigma )$ (and also $\psi$) tending to zero as $t \to \infty$. Performing the change of variables (1.22)
in the auxiliary system (1.20) yields the following modified auxiliary system for the new variables $(q,
\sigma )$
$$\left \{ \begin{array}{l} \partial_t q = i (2t^2)^{-1} \Delta q + t^{-2} (Q(s, q) + Q(\sigma,
W)) + it^{-1} B_{0S} \ q  \\
\\ 
+ i t^{-1} B_S(w, w) q + it^{-1} \left ( 2 B_S(W,q) + B_S(q,q)
\right ) W - R_1(W,S)\\
\\  
\partial_t \sigma = t^{-2} ( s \cdot \nabla \sigma + \sigma \cdot \nabla S) - t^{-1}\nabla
\left ( 2 B_L(W,q) + B_L(q,q)\right ) - R_2 (W, S) \ ,
 \end{array}
\right . \eqno(1.23)$$

\noi where the remainders $R_1(W, S)$ and $R_2(W,S)$ are defined by 
$$R_1(W,S) = \partial_t W - i(2t^2)^{-1} \Delta W - t^{-2} Q(S,W) - it^{-1} (B_{0S} + B_S
(W, W)) W \eqno(1.24)$$
$$R_2(W,S) = \partial_tS - t^{-2}S \cdot \nabla S + t^{-1} \nabla B_{0L}+ t^{-1} \nabla B_L(W,W)
\eqno(1.25)$$

\noi and the dependence of the remainders on $B_0$ has been omitted in the notation. For technical
reasons, it is useful to consider also a partly linearized version of the system (1.23), namely
$$\left \{ \begin{array}{l} \partial_t q' = i (2t^2)^{-1} \Delta q' + t^{-2} (Q(s, q') + Q(\sigma,
W)) + it^{-1} B_{0S} \ q'  \\
\\ 
+ i t^{-1} B_S(w, w) q' + it^{-1} \left ( 2 B_S(W,q) + B_S(q,q)
\right ) W - R_1(W,S)\\
\\  
\partial_t \sigma ' = t^{-2} ( s \cdot \nabla \sigma ' + \sigma \cdot \nabla S) - t^{-1}\nabla
\left ( 2 B_L(W,q) + B_L(q,q)\right ) - R_2 (W, S) \ . 
 \end{array}
\right . \eqno(1.26)$$

The construction of solutions $(w,s)$ of the auxiliary system (1.20) defined for large time and with
prescribed asymptotic behaviour $(W,S)$ proceeds in two steps. The first step consists in solving the
system (1.23) for $(q, \sigma )$ tending to zero at infinity under suitable boundedness properties of
$B_0$ and $(W,S)$ and suitable time decay properties of the remainders $R_1(W,S)$ and $R_2(W,S)$, by a
minor variation of the method used in I. That method consists in first solving the linearized system
(1.26) for $(q' , \sigma ')$ with given $(q , \sigma )$, and then showing that the map $(q, \sigma ) \to
(q', \sigma ')$ thereby defined has a fixed point, by the use of a contraction method. The second step
consists in constructing $(W, S)$ with $W(t)$ tending to $w_+$ as $t \to \infty$ and satisfying the
required boundedness and decay properties. This is done by solving the auxiliary system (1.20) by
iteration to second order as in I and then adding to $W$ an additional term of the same form as that
used in \cite{8r}. The detailed form of $(W, S)$ thereby obtained is too complicated to be given here
and will be given in Section 3 below (see (3.25)-(3.29) and (3.31)).\par

Once the system (1.20) is solved for $(w, s)$, one can proceed therefrom to the construction of a
solution $(u, A)$ of the original WS system. One first defines the phases $\varphi$ and $\phi$ such that
$s = \nabla \varphi$ and $S = \nabla \phi$ and one reconstructs $(u, A)$ from $(w, \varphi )$ by (1.11)
(1.3) (1.5), thereby obtaining a solution of the WS system defined for large time and with prescribed
asymptotic behaviour. The modified wave operator for the WS system is then defined as the map $\Omega :
(u_+, A_+, \dot{A}_+) \to (u, A)$.\par

The main result of this paper is the construction of $(u, A)$ from $(u_+, A_+, \dot{A}_+)$ as described
above, together with the asymptotic properties of $(u, A)$ that follow from that construction. It will
be stated below in full mathematical detail in Proposition 4.1. We give here a heuristic preview of that
result, stripped from most technicalities. We set $\beta = \beta_0 = 1/3$ for definiteness. \\

\noi {\bf Proposition 1.1.} {\it Let $\beta_0 = \beta = 1/3$. Let $(u_+, A_+, \dot{A}_+)$ be such that
$w_+ = Fu_+ \in H^{k_+}$ for sufficiently large $k_+$, that $(A_+, \dot{A}_+)$ be sufficiently regular,
and that $(FA_+, F\dot{A}_+)$ be sufficiently small near $\xi = 0$. Let $(W, S)$ be the approximate
solution of the system (1.20) defined by (3.25)-(3.29) (3.31). Then \par
(1) There exists $T = T(u_+, A_+, \dot{A}_+)$, $1 \leq T < \infty$, such that the auxiliary system
(1.20) has a unique solution ($w, s)$ in a suitable space, defined for $t \geq T$ and such that $(w -W,
s - S)$ tends to zero in suitable norms when $t \to \infty$. \par

(2) There exists $\varphi$ and $\phi$ such that $s = \nabla \varphi$, $S = \nabla \phi$, $\phi (1) = 0$
and such that $\varphi - \phi$ tends to zero in suitable norms when $t \to \infty$. Define $(u, A)$ by
(1.11) (1.3) (1.5). Then $(u, A)$ solves the system (1.1) (1.2) for $t \geq T$ and $(u, A)$ behaves
asymptotically as $(MD \exp (-i \phi ) W$, $A_0 + A_1 (|DW|^2))$ in the sense that the difference tends
to zero in suitable norms (for which each term separately is $O(1)$) when $t \to \infty$.} \\

The unspecified condition that $(FA_+, F\dot{A}_+)$ be sufficiently small near $\xi = 0$ can be shown
to follow from more intuitive conditions in $x$-space, consisting of decay conditions at infinity in
space, and, depending on the values of the parameters defining the relevant function spaces, of some
moment conditions on $(A_+, \dot{A}_+)$.\par

This paper relies on a large amount of material from I. In order to bring out the structure
while keeping duplication to a minimum, we give without proof a shortened logically self-sufficient
sequence of those intermediate results from I that are needed, and we provide a full exposition only
for the parts that are new as compared with I. When quoting I, we shall use the notation (I.p.q) for
equation (p.q) of I and Item I.p.q for Item p.q of I, such as Lemma, Proposition, etc.\par

The remaining part of this paper is organized as follows. In Section 2 we collect notation and some
estimates of a general nature. In Section 3, we study the Cauchy problem at infinity for the
auxiliary system (1.20). We recall from I the existence results of solutions under suitable
boundedness properties of $(W, S)$ and suitable decay properties of the remainders, with the
appropriate modifications (Proposition 3.1). We then define $(W,S)$ and prove that they satisfy the
previous properties, concentrating on the terms in the remainders that are new as compared with I
(Proposition 3.2). We then discuss the assumptions on $(FA_+, F\dot{A}_+)$ at $\xi = 0$ mentioned
above. Finally in Section 4, we construct the wave operators for the WS system (1.1) (1.2) and we
derive the asymptotic properties of the solution $(u, A)$ in their range that follow from the
previous estimates (Proposition 4.1).

\mysection{Notation and preliminary estimates} 
\hspace*{\parindent}
In this section we introduce some notation and we collect a number of estimates which will be used
throughout this paper. We denote by $\parallel \cdot \parallel_r$ the norm in $L^r \equiv
L^r({I\hskip-1truemm R}^{3})$ and we define $\delta (r) = 3/2 - 3/r$.  For
any interval $I$ and any Banach space $X$ we denote by ${\cal C}(I, X)$ the space of strongly continuous
functions from $I$ to $X$ and by $L^{\infty} (I, X)$ the space of measurable essentially bounded
functions from $I$ to $X$. For real numbers $a$ and $b$ we use the notation $a \vee b = {\rm Max}(a,b)$
and $a\wedge b = {\rm Min} (a,b)$. In the estimates of solutions of the relevant equations we shall use
the letter C to denote constants, possibly different from an estimate to the next, depending on various
parameters but not on the solutions themselves or on their initial data. We shall use the notation
$C(a_1, a_2, \cdots )$ for estimating functions, also possibly different from an estimate to the next,
depending on suitable norms $a_1$, $a_2, \cdots$ of the solutions or of their initial data. \par

We shall use the Sobolev spaces $\dot{H}_r^k$ and $H_r^k$ defined for $- \infty < k < + \infty$, $1 \leq
r \leq \infty$ by
$$\dot{H}_r^k = \left \{ u:\parallel u;\dot{H}_r^k\parallel \ \equiv \ \parallel \omega^ku\parallel_r \ <
\infty \right \}$$

\noi and

$$H_r^k = \left \{ u:\parallel u;H_r^k\parallel \ \equiv \ \parallel <\omega>^ku\parallel_r \ <
\infty \right \}$$

\noi where $\omega = (- \Delta)^{1/2}$ and $< \cdot > = (1 + |\cdot |^2)^{1/2}$. The subscript $r$ will
be omitted if $r = 2$.

We shall look for solutions of the auxiliary system (1.20) in spaces of the type ${\cal
C}(I, X^{k, \ell})$ where $I$ is an interval and 
$$X^{k, \ell} = H^k \oplus \omega^{-1} \ H^{\ell}$$
\noi namely
\beq
\label{2.1e}
X^{k, \ell} = \left \{ (w, s) : w \in H^k \ , \ \nabla s \in H^{\ell} \right \} 
\eeq
\noi where it is understood that $\nabla s \in L^2$ includes the fact that $s \in L^6$, and we shall use
the notation 
\beq  
\label{2.2e}
\parallel w; H^k\parallel\ = |w|_k  \ .  \eeq

We shall use extensively the following Sobolev inequalities, stated here in ${I\hskip-1truemm
R}^n$, but to be used only for $n = 3$. \\

\noi {\bf Lemma 2.1.} {\it Let $1 < q$, $r < \infty$, $1 < p \leq \infty$ and $0 \leq j < k$.
If $p = \infty$, assume that $k - j > n/r$. Let $\sigma$ satisfy $j/k \leq \sigma \leq 1$ and
$$n/p - j = (1 - \sigma )n/q + \sigma (n/r - k) \ .$$
\noi Then the following inequality holds
\beq  
\label{2.3e}
\parallel \omega^j u \parallel_p \ \leq C \parallel u \parallel_q^{1 - \sigma} \ \parallel
\omega^ku \parallel_r^{\sigma} \ .  \eeq} 

\indent The proof follows from the Hardy-Littlewood-Sobolev (HLS) inequality (\cite{2r},
p.~117) (from the Young inequality if $p = \infty$), from Paley-Littlewood theory and
interpolation.\\

We shall also use extensively the following Leibnitz and commutator estimates.\\

\noi {\bf Lemma 2.2.} {\it Let $1 < r, r_1, r_3 < \infty$ and
$$1/r = 1/r_1 + 1/r_2 = 1/r_3 + 1/r_4 \ .$$
\noi Then the following estimates hold
\beq  
\label{2.4e}
\parallel \omega^m (uv) \parallel_r \ \leq C \left ( \parallel  \omega^m u \parallel_{r_1} 
 \ \parallel v \parallel_{r_2} + \parallel  \omega^m v \parallel_{r_3} \  
 \parallel u \parallel_{r_4} \right ) \eeq
\noi for $m \geq 0$, and
\beq  
\label{2.5e}
\parallel [\omega^m , u]v\parallel_r \ \leq C \left ( \parallel  \omega^m u \parallel_{r_1} \  
 \parallel v \parallel_{r_2} + \parallel  \omega^{m -1} v \parallel_{r_3} \  
 \parallel \nabla u \parallel_{r_4} \right ) \eeq
\noi for $m \geq 1$, where $[\ ,\ ]$ denotes the commutator.}\\

The proof of those estimates is given in \cite{3r} \cite{4r} with $\omega$ replaced by
$<\omega >$ and follows therefrom by a scaling argument. \\

We next give some estimates of $B_{0L}$, $B_{0S}$, $B_S$ and $B_L$ defined by (1.16) (1.17). It follows
immediately from (1.16) (1.17) that 
\beq  
\label{2.6e}
\parallel \omega^m B_{0L}\parallel_2 \ \leq \left (2 t_0^{\beta_0}\right )^{m-p}  \parallel  \omega^p
B_{0L} \parallel_2 \  \leq \left (2 t^{\beta_0}\right )^{m-p} \parallel \omega^{p} B_0 \parallel_{2} \eeq

\noi for $m \geq p$ and
\beq  
\label{2.7e}
\parallel \omega^m B_{0S}\parallel_2 \ \leq t^{\beta_0(m-p)}  \parallel  \omega^p B_{0S}
\parallel_2 \  \leq t^{\beta_0(m-p)} \parallel \omega^{p} B_0 \parallel_{2} \eeq

\noi for $m \leq p$. Similar estimates hold for $B_L$, $B_S$ with $\beta_0$ replaced by $\beta$. On the
other hand it follows from (1.15) that 
\beq  
\label{2.8e}
\parallel \omega^{m+1} B_{1}(w_1, w_2)\parallel_2 \ \leq I_m \left ( \parallel \omega^{m} (\bar{w}_1
w_2)\parallel_2 \right ) 
\eeq

\noi where $I_m$ is defined by
\beq  
\label{2.9e}
\left ( I_m(f) \right ) (t) = \int_1^{\infty} d\nu \ \nu^{-m-3/2} f(\nu t) \ .
\eeq

We finally collect some estimates of the solutions of the free wave equation\break \noindent $\sq A_0 = 0$
with initial data $(A_+, \dot{A}_+)$ at time zero, given by (1.4). \\

\noi {\bf Lemma 2.3.} {\it Let $k \geq 0$. Let $A_+$ and $\dot{A}_+$ satisfy the conditions
\beq
\label{2.10e}
A_+ , \omega^{-1} \dot{A}_+ \in H^k \qquad , \qquad \nabla^2 A_+, \nabla \dot{A}_+ \in H_1^k \ .
\eeq
Then the following estimate holds~: 
\beq  
\label{2.11e}
\parallel \omega^m A_0 \parallel_r \ \leq \ b_0\ t^{-1+2/r} \qquad {\it for}\  2 \leq r \leq \infty \ ,
\eeq 

\noi for $0 \leq m \leq k$ and for all $t > 0$, where $b_0$ depends on $(A_+, \dot{A}_+)$ through the
norms associated with (\ref{2.10e}).} \\

The estimate (\ref{2.11e}) can be expressed in an equivalent form in terms of $B_0$ defined by
(1.13), namely
\beq  
\label{2.12e}
\parallel \omega^m B_0 \parallel_r \ \leq \ b_0\ t^{m-1/r} \qquad {\it for}\  2 \leq r \leq \infty \ .
\eeq

\noi Furthermore, it follows from (1.17) and (\ref{2.12e}) that
\beq  
\label{2.13e}
\parallel \omega^m B_{0L}\parallel_r \ \leq \parallel  F^{-1} \chi \parallel_1 \ \parallel\omega^m B_{0}
 \parallel_r\  \leq C b_0 \ t^{m-1/r} \eeq

\noi where we have used the Young inequality and the fact that the $L^1$-norm of $F^{-1}\chi$ is 
invariant under the rescaling of $\xi$ by $t^{\beta_0}$ which occurs in (1.17). From (\ref{2.12e})
(\ref{2.13e}) and (1.16) it follows that also  
\beq  
\label{2.14e}
\parallel \omega^m B_{0S}\parallel_r \ \leq  C b_0 \ t^{m-1/r} \ .
\eeq

\noi In the applications, the estimate (\ref{2.12e}) will be used mostly through its consequence
(\ref{2.14e}).

\newpage
\mysection{Cauchy problem at infinity for the auxiliary system}
\hspace*{\parindent}
In this section, we solve the Cauchy problem at infinity for the auxiliary system (1.20) in the
difference form (1.23). We first solve the system (1.23) for $(q, \sigma )$ tending to zero at infinity
under suitable boundedness properties of $(B_0, W, S)$ and suitable time decay properties of the
remainders $R_1(W,S)$ and $R_2(W,S)$. We then construct $(W,S)$ with $W(t)$ tending to $w_+ = Fu_+$ as
$t \to \infty$ and satisfying the required boundedness and decay properties. The method closely follows
that of Sections 6 and 7 of I. \par

We first estimate a single solution of the linearized auxiliary system (1.26) at the level of regularity
where we shall eventually solve the auxiliary system (1.20). The following lemma is basically Lemma
I.6.1, restricted to the case where $1 < k < 3/2$, and sharpened in order to take into account the fact
that the $W$ used in this paper is less regular than that used in I (compare (\ref{3.1e}) below with
(I.6.1)).\\

\noi {\bf Lemma 3.1.} {\it Let $1 < k < \ell < 3/2$ and $\beta >0$. Let $T \geq 1$ and $I = [T, \infty
)$. Let $B_0$ satisfy the
estimate (\ref{2.12e}) for $0 \leq m \leq k$. Let $(U(1/t))W,S) \in {\cal C} (I, X^{k+1,
\ell + 1}) \cap {\cal C}^1(I, X^{k,\ell})$ and let $W$ satisfy
\beq
\label{3.1e}
\mathrel{\mathop {\rm Sup}_{t \in I}} \left \{ \parallel W \parallel_{\infty} \ \vee \parallel W;
H^{3/2}\parallel \ \vee \ t^{1/2-k}  \parallel W; \dot{H}^{k+1}\parallel \right \} \leq a < \infty \ .
\eeq

\noi Let $(q, \sigma )$, $(q', \sigma ') \in {\cal C} (I, X^{k,\ell})$ with $q \in L^{\infty}(I,
H^k) \cap L^2(I,L^2)$ and let $(q', \sigma ')$ be a solution of the system
(1.26) in $I$. Then the following estimates hold for all $t \in I$~:
\bea
\label{3.2e}
\left | \partial_t \parallel q'\parallel_2 \right | &\leq& C \Big \{ t^{-2} a \parallel \nabla
\sigma \parallel_2 + t^{-1-\beta} \ a^2 \ I_0\left ( \parallel q \parallel_2\right ) \nn \\
&&+ t^{-1}\ a \ I_{-1}\left ( \parallel q\parallel_2\ \parallel q\parallel_3 \right ) \Big \} \ +
\ \parallel R_1 (W, S) \parallel_2  \ , \eea
\bea
\label{3.3e}
&&\left | \partial_t \parallel \omega^k q'\parallel_2 \right | \leq C \Big \{ b_0 
\left ( \parallel \omega^{k-1}q'\parallel_2 \ + t^{k-1-\delta /3} \ \parallel q'\parallel_r \right ) \nn
\\ &&+ t^{-2}\ a \left ( \parallel \omega^k \nabla \sigma \parallel_2 \ + t^{k-1/2}\parallel \sigma
\parallel_{\infty} \right ) \nn \\
&&+ t^{-2} \left ( \parallel \nabla s \parallel_{\infty} \ + \  \parallel \omega^{3/2} \nabla s
\parallel_2 \right ) \parallel \omega^{k}q'\parallel_2  \nn \\
&&+ t^{-1} \ a^2 \left ( I_{k-1} \left ( \parallel \omega^{k-1} q \parallel_2 \right )  + \parallel
\omega^{1/2}q'\parallel_2 \right ) \nn \\ 
&&+ t^{-1} \ a \left ( I_{k-1}
\left ( \parallel \omega^{k} q \parallel_2  \ \parallel q \parallel_3 \right ) + I_{1/2} \left ( \parallel
 \omega^{1/2}q\parallel_2 \right ) \parallel \omega^k q'\parallel_2 \right )\nn\\ 
&& + t^{-1} I_{1/2} \left ( \parallel \nabla q \parallel_2^2 \right ) 
\parallel \omega^k q' \parallel_2 \Big \} \ +\ \parallel \omega^k R_1(W,S) \parallel_2  
 \eea

\noi where $s = S + \sigma$ and $0 < \delta = \delta (r) \leq k$. 
\bea
\label{3.4e}
&&\left | \partial_t \parallel \omega^m \nabla \sigma '\parallel_2 \right | \leq C \ t^{-2} \Big
\{ \parallel \nabla s\parallel_{\infty} \ \parallel \omega^m \nabla \sigma ' \parallel_2 \
+ \ \parallel \omega^m \nabla s \parallel_2 \ \parallel \nabla \sigma ' \parallel_{\infty}\nn \\
&&+\ \parallel \omega^m \nabla \sigma\parallel_{2} \ \parallel \nabla S
\parallel_{\infty} \ + \ \parallel \sigma \parallel_{\infty} \ \parallel \omega^m \nabla^2
S\parallel_{2} \Big \}\nn \\
&&+ C \Big \{ t^{-1+ \beta (m+1)} a\ I_0\left ( \parallel q\parallel_2 \right ) + t^{-1 + \beta
(m + 5/2)} \ I_{-3/2} \left ( \parallel q \parallel_2^2 \right ) \Big \} \nn \\
&&+ \ \parallel \omega^m \nabla R_2(W,S)\parallel_2
\eea
\noi for $0 \leq m \leq \ell$}, 
$$\left | \partial_t \parallel \nabla \sigma '\parallel_2 \right | \leq C \ t^{-2} \Big
\{ \parallel \nabla s\parallel_{\infty} \ \parallel \nabla \sigma ' \parallel_2 \
+ \ \parallel \nabla \sigma \parallel_2 \left ( \parallel \nabla S \parallel_{\infty}\ +
\parallel \omega^{3/2} \nabla S \parallel_2 \right ) \Big \} $$
$$+ \ C \Big \{ t^{-1+\beta}\ a\ I_0 \left ( \parallel q \parallel_2\right ) + t^{-1+5\beta /2}\
I_{-3/2} \left ( \parallel q \parallel_2^2 \right ) \Big \} \ + \ \parallel \nabla
R_2(W,S)\parallel_2 \ . \eqno(3.4)_0$$
\par \vskip 5 truemm

\noi {\bf Proof}. (\ref{3.2e}) is identical with (I.6.2) and is proved in the same way. In order to prove
(\ref{3.3e}), we start from (cf. I.6.9))

\bea
\label{3.5e}
&&\left | \partial_t \parallel \omega^k q'\parallel_2 \right | \leq t^{-1}
\parallel [\omega^{k}, B_{0S}]q'\parallel_2 \ +\  t^{-2} \left (  \parallel
[\omega^k,s]\cdot \nabla q'\parallel_2 \right .\nn \\ 
&&\left . + \parallel (\nabla \cdot s)\omega^k q' \parallel_2 \ + \ \parallel \omega^k \left (
(\nabla \cdot s) q'\right ) \parallel_{2} \ + \ \parallel \omega^k Q(\sigma, W) \parallel_2
\right ) \nn \\
 && t^{-1} \left ( \parallel [\omega^k ,B_S(w,w)]q' \parallel_{2} \ + \ \parallel
\omega^k \left ( 2B_S(W,q) + B_S(q,q)\right ) W \parallel_2 \right )  \nn \\
&&+\ \parallel \omega^k \ R_1(W,S) \parallel_2 
 \eea

\noi and we estimate the various terms in the RHS successively. \par

The contribution of $B_0$ is estimated exactly as in I and yields 
\beq
\label{3.6e}
\parallel [\omega^k, B_{0S}]q'\parallel_2\ \leq \ C\ b_0 \left ( t \parallel \omega^{k-1}
q'\parallel_2 \ + t^{k-\delta /3} \ \parallel q'\parallel_r \right ) \ .\eeq
\indent The contribution of $Q(s, q')$ is estimated by Lemmas 2.1 and 2.2 as 
\bea
\label{3.7e}
&&\parallel [\omega^k, s]\cdot \nabla q'\parallel_2\ + \ \parallel
(\nabla \cdot s) \omega^k q'\parallel_2 \ + \ \parallel \omega^k \left ( (\nabla \cdot
s)q'\right ) \parallel_2 \nn \\
&&\leq C \left ( \parallel \nabla s\parallel_{\infty} \ + \ \parallel
\omega^{3/2} \nabla s \parallel_2 \right ) \parallel \omega^k q' \parallel_2 
\eea
\noi in the same way as in I, in the case $k < 3/2$. \par

The contribution of $Q(\sigma , W)$ is estimated by Lemma 2.1 and 2.2 as 
$$\parallel \omega^kQ(\sigma , W)\parallel_2\ \leq C \left (  \parallel\sigma \parallel_{\infty} 
\ \parallel \omega^k \nabla W\parallel_2 \ + \  \parallel \omega^k \sigma \parallel_6\ \parallel
\nabla W\parallel_3 \right .$$
$$\left . + \ \parallel \omega^k\nabla \sigma \parallel_{2} \ \parallel
W\parallel_{\infty} \ + \ \parallel \nabla \sigma \parallel_r \ \parallel
\omega^kW \parallel_{3/k} \right )$$

\noi with $\delta (r) = k$,
\bea
\label{3.8e}
&&\cdots \leq C \left (  \parallel \sigma \parallel_{\infty} \ \parallel \omega^k \nabla  W 
\parallel_2 \ + \ \parallel \omega^k \nabla \sigma \parallel_{2} \left ( \parallel W \parallel_{\infty}
\ + \ \parallel\omega^k W \parallel_{3/k} \right ) \right ) \nn \\ 
&&\leq C\ a \left ( t^{k-1/2} \parallel
\sigma \parallel_{\infty} \ + \ \parallel \omega^k \nabla \sigma \parallel_2 \right )
\eea
\noi by (\ref{3.1e}).\par

The contribution of $B_S$ with $w = W + q$ yields a number of terms which we order by
increasing powers of $q$, $q'$. We first expand
\beq
\label{3.9e}
B_S(w, w) = B_S(W,W) + 2B_S(W,q) + B_S(q,q) \ .
\eeq
\noi By Lemmas 2.1 and 2.2, we estimate
$$\parallel [\omega^k, B_S(W,W)]q'\parallel_2\ \leq C \left (  \parallel\nabla B_S(W,W)
\parallel_{3/\varepsilon} + \ \parallel \omega^{k} B_S(W,W) \parallel_{3/(k-1+
\varepsilon )}\right )$$
$$\times \parallel \omega^{k-1+ \varepsilon} q'\parallel_2$$

\noi for $\varepsilon > 0$. Taking $\varepsilon = 3/2 - k$ yields
\bea
\label{3.10e}
&&\cdots \leq C \parallel \omega^{k+1} B_1(W,W)  \parallel_2 \ \parallel \omega^{1/2} q' \parallel_{2} \nn
\\ 
&&\leq C\ I_k \left ( \parallel \omega^kW \parallel_{2} \ \parallel W \parallel_{\infty}
\right ) \parallel \omega^{1/2} q' \parallel_{2}\nn \\
&&\leq C\ a^2 \parallel \omega^{1/2} q' \parallel_2 
\eea

\noi by Lemma 2.2 again and by (\ref{2.8e}) (\ref{3.1e}).\par

In a similar way, we estimate by Lemmas 2.1, 2.2 and by (\ref{2.8e})
$$\parallel [\omega^k, B_S(W,q)]q'\parallel_2\ \leq C \left (  \parallel\nabla B_S(W,q)
\parallel_{3} \ \parallel \omega^{k-1} q' \parallel_{6} \ + \ \parallel \omega^{k} B_S(W,q)
\parallel_{3/k} \parallel q'\parallel_r \right ) $$

\noi with $\delta (r) = k$, 
\bea
\label{3.11e}
&&\cdots \leq C \parallel \omega^{3/2} B_1(W,q)  \parallel_2 \ \parallel \omega^{k} q' \parallel_{2} \nn
\\ 
&&\leq C\ I_{1/2} \left ( \parallel W \parallel_{\infty} \ \parallel \omega^{1/2}q \parallel_{2} \ +\
\parallel \omega^{1/2} W \parallel_{6} \parallel q \parallel_{3} \right ) \parallel \omega^{k} q'
\parallel_{2}\nn \\ &&\leq C\ a\ I_{1/2} \left ( \parallel \omega^{1/2} q \parallel_2 \right ) \parallel
\omega^{k} q' \parallel_{2} \ .  
\eea

\noi In a similar way, we estimate
\bea
\label{3.12e}
&&\parallel [\omega^k, B_S(q,q)]q'\parallel_2\ \leq C \left (  \parallel\nabla B_S(q,q)
\parallel_{3} \ + \ \parallel \omega^{k} B_S(q,q)\parallel_{3/k} \right ) \parallel \omega^{k} q'
\parallel_{2} \nn \\
&&\leq C \parallel \omega^{3/2} B_1(q,q)  \parallel_2 \ \parallel \omega^{k} q' \parallel_{2} \nn
\\
&&\leq C\ I_{1/2} \left ( \parallel \omega^{1/2} q \parallel_{3} \ \parallel q
\parallel_{6} \right )  \parallel \omega^{k} q' \parallel_{2} \nn \\
&&\leq C\ I_{1/2} \left ( \parallel \nabla q \parallel_2^2 \right ) \parallel \omega^{k} q'
\parallel_{2} \ .  
\eea

\noi We next estimate in a similar way
\bea
\label{3.13e}
&&\parallel \omega^k (B_S(W,q)W)\parallel_2\ \leq C \parallel\omega^k B_1(W,q)
\parallel_{2} \left ( \parallel W \parallel_{\infty} + \parallel \omega^{k} W
\parallel_{3/k} \right )\nn \\
&&\leq C \ a \parallel \omega^{k} B_1(W,q)  \parallel_2  \nn
\\
&&\leq C\ a\ I_{k-1} \left ( \left ( \parallel W  \parallel_{\infty} \ + \ \parallel \omega^{k-1} W
\parallel_{3/(k-1)} \right )  \parallel \omega^{k-1} q \parallel_{2} \right )\nn \\
&&\leq C\ a^2 \ I_{k-1} \left ( \parallel \omega^{k-1} q \parallel_2\right )  \ .  
\eea

\noi Finally, we estimate in a similar way
\bea
\label{3.14e}
&&\parallel \omega^k (B_S(q,q)W)\parallel_2\ \leq C \parallel\omega^k B_1(q,q)
\parallel_{2} \left ( \parallel W \parallel_{\infty} + \parallel \omega^{k} W
\parallel_{3/k} \right )\nn \\
&&\leq C\ a\ I_{k-1} \left ( \parallel \omega^{k-1} |q|^2  \parallel_{2} \right ) \nn \\
&&\leq C\ a \ I_{k-1} \left ( \parallel \omega^{k} q \parallel_2 \ \parallel q
\parallel_3\right )  \ .   \eea

\noi Substituting (\ref{3.6e})-(\ref{3.8e}) and (\ref{3.10e})-(\ref{3.14e}) into (\ref{3.5e}) yields
(\ref{3.3e}). \par

The estimates (\ref{3.4e}) and (3.4)$_0$ of $\sigma '$ are identical with (I.6.4) and (I.6.4)$_0$ and
have exactly the same proof. The additional term with $B_{0L}$ in the equation for $\sigma '$ is
included in $R_2(W,S)$ and does not appear explicitly at this stage.\par\nobreak
\hfill $\sq$\par

From there on, the treatment of the Cauchy problem at infinity for the auxiliary system follows that
given in I verbatim. We need to estimate the difference of two solutions of the linearized auxiliary
system (1.26), and that estimate, given by Lemma I.6.2, requires no modification because it uses
regularity properties of $W$ which are weaker than (\ref{3.1e}). We then solve the Cauchy problem first
for the linearized auxiliary system (1.26) with finite initial time by Proposition I.6.1, then at
infinity by Proposition I.6.2, and then for the auxiliary system (1.20) or (1.23) by a contraction
method, by Proposition I.6.3, part (2). The only difference in the proof of Propositions I.6.2 and
I.6.3 is due to the term 
$$t^{-2}a\ t^{k-1/2} \parallel \sigma \parallel_{\infty}$$ 
\noi in (\ref{3.3e}), which did not appear in Lemma I.6.1, and which is due to the fact that the
assumption (\ref{3.1e}) is weaker than (I.6.1). That term generates an additional term
$$a\ Z\ t^{-1-\lambda - 3(1 - \beta)/2}$$
\noi in the RHS of (I.6.59), with time decay strictly better than $t^{-1-\lambda}$ and therefore
harmless. \par

We now state the first main result of this section, which corresponds to Proposition I.6.3, part (2).\\

\noi {\bf Proposition 3.1.} {\it Let $1 < k < 3/2 < \ell$. Let $\lambda_0$, $\lambda$ and $\beta$ satisfy
the conditions
\beq
\label{3.15e}
\lambda > 0 \qquad , \qquad (1<) \lambda + k < \lambda_0 \ ,
\eeq
\beq
\label{3.16e}
0 < \beta < 2/3 \qquad , \qquad \beta (\ell + 1 ) < \lambda_0 \ .
\eeq
\noi Let $(A_+, \dot{A}_+)$ satisfy the conditions (2.10). Let $(U(1/t)W,S) \in {\cal C}([1, \infty),
X^{k+1, \ell + 1}) \cap {\cal C}^1([1, \infty), X^{k,\ell})$ and let $(W, S)$ satisfy the estimates 
\beq
\label{3.17e}
\mathrel{\mathop {\rm Sup}_{t \geq 1}} \left \{ \parallel W \parallel_{\infty} \ \vee \ \parallel W;
H^{3/2}\parallel \ \vee \ t^{1/2 - k}\parallel W; \dot{H}^{k+1} \parallel \right \} \leq a < \infty \ , 
\eeq
\beq
\label{3.18e}
\parallel \omega^m \nabla S \parallel_2 \ \leq b\ t^{1- \eta + \beta (m - 3/2)}
\eeq
\noi for some $\eta > 0$ and for $0 \leq m \leq \ell + 1$,
\beq
\label{3.19e}
\parallel R_1(W,S) \parallel_2 \ \leq c_0 \ t^{-1-\lambda_0} \ ,
\eeq
\beq
\label{3.20e}
\parallel \omega^k R_1(W,S)\parallel_2 \ \leq c_1 \ t^{-1-\lambda} \ ,
\eeq
\beq
\label{3.21e}
\parallel \omega^m \nabla R_2(W,S)\parallel_2 \ \leq c_2 \ t^{-1-\lambda_0 + \beta (m+1)} \ {\it for} \
0 \leq m \leq \ell \ . 
\eeq

\noi Then there exists $T$, $1 \leq T < \infty$ and positive constants $Y_0$, $Y$ and $Z$, depending on
$(A_+, \dot{A}_+)$ through the norms in (\ref{2.10e}) and depending on $k$, $\ell$, $\beta$, $\lambda_0$,
$\lambda$, $a$, $b$, $c_0$, $c_1$ and $c_2$ such that the auxiliary system (1.20) has a unique solution
$(w, s) \in {\cal C}(I, X^{k,\ell})$, where $I = [T, \infty )$, satisfying the estimates
\beq
\label{3.22e}
\parallel w - W \parallel_2 \ \leq Y_0\ t^{-\lambda_0} \ ,
\eeq
\beq
\label{3.23e}
\parallel \omega^k (w - W) \parallel_2 \ \leq Y\ t^{-\lambda} \ ,
\eeq
\beq
\label{3.24e}
\parallel \omega^m \nabla  (s - S) \parallel_2 \ \leq Z\ t^{-\lambda_0 + \beta (m+1)} \ {\it for} \
0 \leq m \leq \ell \ , 
\eeq
\noi for all $t \in I$.} \\

We now turn to the construction of approximate solutions $(W, S)$ of the system (1.20) satisfying the
assumptions of Proposition 3.1 and in particular the estimates (\ref{3.17e}) (\ref{3.18e}) of $(W, S)$
and the estimates (\ref{3.19e})-(\ref{3.21e}) of the remainders. In I we took for $(W, S)$ the second
order approximate solution of the system (1.20) in an iterative scheme not taking into account the
terms containing $B_0$, thereby ending with an explicit $B_0W$ term in the remainder $R_1(W, S)$.
Here, following \cite{8r}, we improve that asymptotic form by adding one more term in $W$, so as to
partly cancel the $B_{0S}W$ term in $R_1(W,S)$. Thus we define
\beq
\label{3.25e}
W = w_0 + w_1 + w_2 \equiv W_1 + w_2 \quad , \quad S = s_0 + s_1 
\eeq  

\noi where $w_0$, $s_0$, $w_1$, $s_1$ are the same as in I, namely
\beq
\label{3.26e}
w_0 = U^*(1/t) w_+ \ ,
\eeq
\beq
\label{3.27e}
s_0(t) = - \int_1^t dt'\ t'^{-1} \nabla B_L(w_0(t'),w_0(t')) \ ,
\eeq
\beq
\label{3.28e}
w_1(t) = - U^*(1/t) \int_t^{\infty} dt' \ t'^{-2} U(1/t') Q(s_0(t'),w_0(t')) \ ,
\eeq
\beq
\label{3.29e}
s_1(t) = - \int_t^{\infty} dt'\ t'^{-2} s_0(t') \cdot \nabla s_0(t') + 2 \int_t^{\infty} dt'\ t'^{-1}
\nabla B_L (w_0(t'), w_1 (t'))\ .
\eeq

\noi In order to partly cancel $B_{0S}W$ in $R_1(W,S)$, we take $w_2 = hw_0$, thereby obtaining a linear
contribution of $w_2$ to $R_1(W,S)$
\bea
\label{3.30e}
\left ( \partial_t - i(2t^2)^{-1} \Delta \right ) w_2 &=& h\left ( \partial_t - i(2t^2)^{-1} \Delta
\right ) w_0 - i t^{-2} \nabla h \cdot \nabla w_0 + (\partial_t h) w_0\nn \\
 &&- i ( 2t^2)^{-1} (\Delta
h)w_0 \ . \eea

\noi The first term in the RHS is small, actually zero, by the choice of $w_0$. We use the last term in
the RHS to cancel the main contribution $B_{0S}w_0$ of $B_{0S}W$ by making the choice
\beq
\label{3.31e}
w_2 = h w_0 \qquad , \qquad h = - 2 t \Delta^{-1} B_{0S}\ .
\eeq

\noi Note that because of the short range cut-off in $B_{0S}$, $h$ is well defined, actually $h \in
{\cal C}([1, \infty), H^{k+2})$. With that choice, the remainders become
\beq
\label{3.32e}
R_i(W,S) = R_{i0}(W,S) + R_{i\nu}(W,S) \qquad i = 1,2 \ ,
\eeq

\noi where $R_{i0}(W,s)$ are the parts not containing $w_2$ or $B_{0L}$, namely
\bea
\label{3.33e}
R_{10}(W,S) &=& - t^{-2} \left ( Q(S,w_1) + Q (s_1, w_0) \right ) \nn \\
&&- i t^{-1} \left ( B_{0S}w_1 + B_S (W_1, W_1) W_1 \right )
\eea 
\beq
\label{3.34e}
R_{20}(W,S) = - t^{-2} \left ( s_0 \cdot \nabla s_1 + s_1 \cdot \nabla s_0 + s_1 \cdot \nabla s_1 \right
) + t^{-1} \nabla B_L (w_1,w_1)\ ,
\eeq

\noi while $R_{i\nu}(W,S)$ are the parts containing $w_2$ or $B_{0L}$, namely
\bea
\label{3.35e}
R_{1\nu}(W,S) &=& - t^{-2} Q(S,w_2) - i t^{-1} B_{0S} w_2 + (\partial_t h) w_0 - i t^{-2} \nabla h
\cdot \nabla w_0 \nn \\ 
&&- i t^{-1} \left ( B_{S}(W,W)w_2 + B_S (W + W_1, w_2) W_1 \right
) \ , \eea
\beq
\label{3.36e}
 R_{2\nu}(W,S) = t^{-1} \nabla B_{0L} + t^{-1} \nabla B_L(W+W_1, w_2) \ .
\eeq

\noi The parts $R_{i0}$ of the remainders are the remainders occurring in I, up to the replacement of
$B_0$ by $B_{0S}$ and the disappearance of the term $B_{0S}w_0$, precisely the term which was
responsible for the support condition in I. Up to a minor point (see below), $(W_1, S)$ and
$R_{i0}(W,S)$ have been estimated in I as follows (see Lemma I.7.1).\\

\noi {\bf Lemma 3.2.} {\it Let $0 < \beta < 1$, $k_+ \geq 3$, $w_+ \in H^{k_+}$ and $a_+ =
|w_+|_{k_+}$. Then the following estimates hold for all $t \geq 1$~:
\beq
\label{3.37e}
|w_0|_{k_+} \ \leq a_+ \ ,
\eeq

\beq
\label{3.38e}
\parallel \omega^m \ s_0\parallel_2 \ \leq \left \{ \begin{array}{ll} C\ a_+^2 \
\ell n \ t &\qquad \hbox{for}\ 0 \leq m \leq k_+ \\
\\
C \ a_+^2 \ t^{\beta (m - k_+)} &\qquad \hbox{for}\ m > k_+ \ ,
\end{array} \right .  
\eeq

\beq
\label{3.39e}
|w_1|_{k_+ - 1} \leq C \ a_+^3\ t^{-1} (1 + \ell n\ t) \ ,
\eeq

\beq
\label{3.40e}
\parallel \omega^m \ s_1\parallel_2 \ \leq \left \{ \begin{array}{l} C\ a_+^4 \ t^{-1} (1
+ \ell n \ t)^2 \quad \qquad \hbox{for}\ 0 \leq m \leq k_+ - 1 \\ \\
C \ a_+^4 \ t^{-1 + \beta (m +1- k_+)} (1 + \ell n \ t) \\
\qquad \hbox{for}\ k_+ - 1 < m < k_+ -1 +
\beta^{-1} \ , \end{array} \right . 
\eeq

\beq
\label{3.41e}
\parallel \omega^m \ R_{20}(W,S)\parallel_2 \ \leq \left \{ \begin{array}{l} C(a_+) \ t^{-3} (1
+ \ell n \ t)^3 \quad \qquad \hbox{for}\ 0 \leq m \leq k_+ - 2 \\ \\
C(a_+) \ t^{-3 + \beta (m +2- k_+)} (1 + \ell n \ t)^2 \\ \qquad \hbox{for}\ k_+ - 2 < m < k_+
-2 + \beta^{-1} \ , \end{array} \right .
\eeq

\noi Let in addition $0 \leq k \leq k_+ - 1$ and let $B_0$ satisfy the estimate (\ref{2.12e}) for $0
\leq m \leq k$. Then}
\bea
\label{3.42e}
\parallel \omega^m R_{10}(W,S)\parallel_2  &\leq& C(a_+)
\left \{ t^{-3} (1 + \ell n \ t)^2 + t^{-1-\beta (k_+ - m + 1)} \right \} \nn \\
&&+ C\ b_0 \ a_+^3 \ t^{m-5/2}(1 + \ell n\ t) \qquad \hbox{for}\ 0 \leq m \leq
k \ .
\eea
\vskip 5 truemm
\noi {\bf Proof.} The estimates (\ref{3.37e})-(\ref{3.42e}) are those of Lemma I.7.1 except for the
estimate of the term $t^{-1} B_{0S} w_1$ in $R_{10}(W,S)$ which is responsible for the last term in
(\ref{3.42e}). We estimate that term by Lemma 2.1 and 2.2 and by (\ref{2.12e}) (\ref{3.39e}) as 
\bea
\label{3.43}
\parallel \omega^m B_{0S}w_1\parallel_2 &\leq& \left \{ \begin{array}{ll} C\left ( \parallel \omega^m
B_{0}\parallel_2\  \parallel w_1  \parallel_{\infty}\ + \  \parallel B_{0S}  \parallel_r \ 
\parallel\omega^m w_1 \parallel_{1/m} \right ) &\hbox{for}\ m \leq 1/2 \\ & \\ C\left ( \parallel \omega^m
B_{0}\parallel_2\  \parallel w_1  \parallel_{\infty}\ + \  \parallel B_{0S}  \parallel_{\infty} \ 
\parallel\omega^m w_1 \parallel_{2} \right ) &\hbox{for}\ m > 1/2 \end{array} \right . \nn \\
&& \nn \\
&\leq& C\ b_0 \  a_+^3 \ t^{m -3/2} (1 + \ell n \ t) 
\eea

\noi for $0 \leq m \leq k$, and $1/r = 1/2 - m$ for $m \leq 1/2$, which completes the proof of
(\ref{3.42e}). \par\nobreak
\hfill $\sq$\par

We now turn to estimating $R_{i\nu}(W,S)$, $i = 1,2$. We first reduce that question to that of
estimating $h$ and $B_{0L}$, assuming for the moment a boundedness property of $w_2$ which is part of
(\ref{3.17e}) and which we shall prove later. We define the auxiliary space
\beq
\label{3.44e}
Y = L^{\infty} \left ( [1, \infty), L^{\infty} \cap H^{3/2} \right )
\eeq

\noi and we remark that for $k_+ > 5/2$, it follows from (\ref{3.37e}) (\ref{3.39e}) that $W_1 = w_0 +
w_1 \in Y$. We can now state the estimates of $R_{i\nu}(W, S)$. \\

\noi {\bf Lemma 3.3.} {\it Let $1 < k < 3/2$ and  $k_+ \geq 3$, let $w_+ \in H^{k_+}$ and let $a_+ =
|w_+|_{k_+}$. Assume that $w_2 \in Y$ and let
\beq
\label{3.45e}
\parallel W_1\parallel_Y \ \vee \ \parallel W \parallel_Y \ \leq a < \infty \ .
\eeq
\noi Let $B_0$ satisfy the estimate (\ref{2.12e}) for $0 \leq m \leq k$. Then the following estimates
hold for all $t \geq 1$~:
\bea
\label{3.46e}
\parallel R_{1\nu}(W,S)\parallel_2 &\leq& C(a_+) t^{-2} \ell n \ t \parallel\nabla h \parallel_2
+ C\ a_+ \Big \{ b_0 \ t^{-1} \parallel h \parallel_2 \nn \\  
&&+ \ \parallel \partial_t h \parallel_2 \ + a^2\  t^{-1} \left ( \parallel h \parallel_2 + I_0
(\parallel h \parallel_2) \right ) \Big \} \ , 
\eea
\bea
\label{3.47e}
&&\parallel \omega^k R_{1\nu}(W,S)\parallel_2 \ \leq C(a_+) t^{-2} \ell n \ t \left ( \parallel \omega^k
\nabla h \parallel_2 \ + \ \parallel \omega^{\delta} \nabla h \parallel_2 \right ) \nn \\ 
&&+ C\ a_+ \Big \{ b_0 \left ( t^{k-1-\delta /3} \parallel \omega^{\delta} h \parallel_2 \ + t^{-1}
\parallel \omega^k h \parallel_2 \right ) + \ \parallel\omega^k \partial_t h \parallel_2 \ + \ \parallel
\omega^{\delta} \partial_t h \parallel_2 \nn \\
&&+ t^{-2} \left ( \parallel \omega^{k+1} h \parallel_2 \ + \ \parallel \omega^2 h \parallel_2 \right )
+ t^{-1} \ a^2 \left ( \parallel \omega^kh\parallel_2 \ + I_{k-1} \left ( \parallel \omega^{k-1} h
\parallel_2 \right ) \right ) \Big \}\nn \\
\eea

\noi for $0 < \delta < 1/2$,
\beq
\label{3.48e}
\parallel \omega^m R_{2\nu}(W,S)\parallel_2 \ \leq t^{-1} \parallel \omega^{m+1} B_{0L}\parallel_2 \ +
C\ a \ a_+\ t^{-1 + \beta m} I_0 \left ( \parallel h \parallel_2 \right )  
\eeq

\noi for all $m \geq 0$.}\\

\noi {\bf Proof.} We first consider $\parallel R_{1\nu}(W,S)\parallel_2$. We estimate successively
\bea
\label{3.49e}
\parallel Q(S, w_2) \parallel_2 &\leq& \parallel S \parallel_{\infty} \left ( \parallel w_0
\parallel_{\infty} \ \parallel \nabla h \parallel_{2} \ + \ \parallel \nabla w_0\parallel_3 \ \parallel
h \parallel_6 \right )\nn \\
&&+\ \parallel \nabla S \parallel_3 \ \parallel w_0 \parallel_{\infty} \ \parallel h \parallel_6 \ \leq
C(a_+) \ell n \ t \parallel\nabla h \parallel_2  \eea

\noi by (\ref{3.37e}) (\ref{3.38e}) (\ref{3.40e}) and Lemma 2.1,
\beq
\label{3.50e}
\parallel B_{0S} w_2 \parallel_2 \ \leq \ \parallel B_{0S} \parallel_{\infty} \ \parallel w_0
\parallel_{\infty} \ \parallel h \parallel_2 \ \leq C\ b_0 \ a_+ \parallel h\parallel_2 \ ,
 \eeq

\noi by (\ref{2.14e}) and (\ref{3.37e}),
\beq
\label{3.51e}
\parallel(\partial_t h)w_0\parallel_2 \ \leq \ \parallel w_0 \parallel_{\infty}\ \parallel \partial_t h
\parallel_2 \ \leq C\ a_+ \parallel\partial_t h \parallel_2 \ ,
 \eeq
\beq
\label{3.52e}
\parallel \nabla h \cdot \nabla w_0\parallel_2\ \leq \ \parallel \nabla w_0 \parallel_{\infty} \
\parallel\nabla h \parallel_2\ \leq C \ a_+ \parallel \nabla h \parallel_2 \ , 
\eeq
\beq
\label{3.53e}
\parallel B_S (W,W)w_2\parallel_2 \ \leq \ \parallel B_S(W,W) \parallel_{\infty} \ \parallel w_0
\parallel_{\infty} \ \parallel h \parallel_2 \ \leq C\ a_+\ a^2 \parallel h \parallel_2 \ , 
\eeq

\noi by estimating $B_S(W, W)$ in a way similar to that in Lemma 3.1,
\bea
\label{3.54e}
&&\parallel B_S (W+ W_1,w_2)W_1\parallel_2 \ \leq C\parallel W_1 \parallel_{3} \ \parallel \omega B_1(W
+ W_1, w_2) \parallel_{2} \nn \\
&&\leq C \parallel W_1\parallel_3 \ \parallel W + W_1 \parallel_{\infty} \ \parallel w_0
\parallel_{\infty} \ I_0(\parallel h \parallel_2) \leq C\ a_+ \ a^2\ I_0 (\parallel h \parallel_2 )  
\eea

\noi by Lemma 2.1 and by (\ref{2.8e}). Collecting (\ref{3.49e})-(\ref{3.54e}) yields (\ref{3.46e}). \par

We next consider $\omega^k R_{1\nu}(W,S)$. We estimate successively
\bea
\label{3.55e}
&&\parallel \omega^k Q(S, w_2) \parallel_2 \ \leq C \Big \{ \parallel w_0 \parallel_{*} \ 
\parallel S \parallel_{*} \ \parallel \omega^{k+1} h \parallel_{2} \ + \Big ( \parallel w_0
\parallel_* \ \parallel \omega^{k+1} S \parallel_{3/(1 + \delta )} \nn \\
&&+\ \parallel \omega^{k+1} w_0  \parallel_{3/(1 + \delta )} \ \parallel S \parallel_{*} \Big )
\parallel \omega^{\delta} \nabla h \parallel_2\Big \} \nn \\
&&  \leq C(a_+) \ell n \ t \left ( \parallel
\omega^k \nabla h \parallel_2 \ + \ \parallel \omega^{\delta} \nabla h \parallel_2 \right ) 
\eea

\noi by (\ref{3.37e}) (\ref{3.38e}) (\ref{3.40e}) and Lemmas 2.1 and 2.2, with
$$\parallel f \parallel_* \ = \ \parallel f \parallel_{\infty}\ + \ \parallel \nabla f \parallel_3 \ ,$$
\bea
\label{3.56e}
&&\parallel \omega^k (B_{0S} w_2) \parallel_2 \ \leq C \Big \{ \parallel B_0 \parallel_{\infty} \ 
\parallel w_0 \parallel_{\infty} \ \parallel \omega^{k} h \parallel_{2} \ + \Big ( \parallel \omega^k B_0
\parallel_{3/\delta} \ \parallel w_0 \parallel_{\infty} \nn \\
&&+\ \parallel B_0  \parallel_{\infty} \ \parallel \omega^k w_0 \parallel_{3/\delta } \Big )
\parallel \omega^{\delta} h \parallel_2\Big \}\nn \\
&&   \leq C\ b_0 \ a_+ \left ( \parallel
\omega^k h \parallel_2 \ + t^{k-\delta /3}\ \parallel \omega^{\delta} h \parallel_2 \right ) 
\eea

\noi by Lemmas 2.1 and 2.2 and by (\ref{2.14e}) (\ref{3.37e}),
\bea
\label{3.57e}
\parallel \omega^k ((\partial_t h ) w_0) \parallel_2 &\leq& C \left ( \parallel w_0 \parallel_{\infty} \
\parallel \omega^k \partial_t h \parallel_2  \ + \ \parallel \omega^k w_0 \parallel_{3/\delta} \
\parallel\omega^{\delta} \partial_t h \parallel_2 \right )  \nn \\
&\leq& C\ a_+ \left ( \parallel \omega^k \partial_t h \parallel_2 \ + \ \parallel \omega^{\delta}
\partial_t h \parallel_2 \right )
\eea

\noi by Lemmas 2.1 and 2.2 and by (\ref{3.37e}),
\bea
\label{3.58e}
\parallel \omega^k (\nabla h \cdot \nabla w_0) \parallel_2 &\leq& C \left ( \parallel \nabla w_0
\parallel_{\infty} \ \parallel \omega^{k+1} h \parallel_2  \ + \ \parallel \omega^{k+1} w_0
\parallel_{3} \ \parallel \nabla h  \parallel_6 \right )  \nn \\
&\leq& C\ a_+ \left ( \parallel \omega^{k+1} h \parallel_2 \ + \ \parallel \omega^{2}
h \parallel_2 \right )
\eea

\noi by Lemmas 2.1 and 2.2 and by (\ref{3.37e}). We next estimate
\bea
\label{3.59e}
&&\parallel \omega^k (B_S(W,W)w_2) \parallel_2 \ \leq C \left ( \parallel B_S w_0
\parallel_{\infty} \ + \ \parallel \omega^{k} (B_S w_0) \parallel_{3/k}  \right ) \parallel \omega^{k} h
\parallel_{2}\nn \\
&&\leq C \left ( \parallel B_S \parallel_{\infty} \ + \ \parallel \omega^{k} B_S \parallel_{3/k}  \right
) \left ( \parallel w_{0} \parallel_{\infty}\ + \ \parallel \omega^{k} w_0
\parallel_{3/k}\right ) \parallel \omega^{k} h \parallel_{2} \nn \\
&&\leq C\ a_+ \ a^2 \parallel \omega^k h \parallel_2
\eea

\noi where we have omitted the arguments in $B_S$, by Lemmas 2.1 and 2.2 and by (\ref{3.37e})
(\ref{3.45e}), and after estimating $B_S$ in a way similar to that in Lemma 3.1. In the same way
\bea
\label{3.60e}
&&\parallel \omega^k (B_S(W +W_1,w_2)W_1) \parallel_2 \ \leq C \left ( \parallel W_1
\parallel_{\infty} \ + \ \parallel \omega^{k} W_1 \parallel_{3/k}  \right ) \parallel \omega^{k} B_1(W +
W_1, w_2) \parallel_{2}\nn \\
&&\leq C \ a\ I_{k-1} \left ( \parallel \omega^{k-1} h \parallel_{2} \left (  \parallel w_0(W + W_1) 
\parallel_{\infty}  \ + \ \parallel \omega^{k-1} (W + W_1)w_0\parallel_{3/(k-1)}\right ) \right )\nn \\
&&\leq C\ a_+ \ a^2 \ I_{k-1} \left ( \parallel \omega^{k-1} h \parallel_2 \right ) \ .
\eea

\noi Collecting (\ref{3.55e})-(\ref{3.60e}) yields (\ref{3.47e}). \par

We finally estimate $R_{2\nu}$. From (\ref{2.6e}) (\ref{2.8e}) we obtain 
$$\parallel \omega^m R_{2\nu}(W, S) \parallel_2 \ \leq t^{-1} \parallel \omega^{m+1} B_{0L}
\parallel_2$$
$$+ C\ t^{-1 + \beta m}\   I_0 \left ( \parallel W + W_1 \parallel_{\infty} \ \parallel w_0
\parallel_{\infty} \ \parallel h \parallel_2 \right )$$

\noi which yields (\ref{3.48e}) by using (\ref{3.37e}) (\ref{3.45e}). \par \nobreak \hfill $\sq$\par

In order to complete the estimate of the parts $R_{i\nu}(W,S)$, $i = 1,2$, of the remainders, we now
estimate $h$ and $B_{0L}$. Those estimates require some restrictions on the behaviour of
$(FA_+, F\dot{A}_+)$ at $\xi = 0$. Those restrictions are imposed in a dilation
homogeneous way through the use of a parameter $\mu \in (-1, 1)$ in terms of quantities which have the
same scaling properties as $\parallel  A_+ ; \dot{H}^{-3/2-\mu}\parallel$ and $\parallel\dot{A}_+ ;
\dot{H}^{-5/2 - \mu} \parallel$. They will be further discussed at the end of this section.\\

\noi {\bf Lemma 3.4.} {\it Let $1 < k < 3/2$ and $-1 < \mu < 1$. Let $(A_+, \dot{A}_+)$ satisfy the
conditions
\beq
\label{3.61e}
A_+ \in H^{k-1} \qquad , \qquad \dot{A}_+ \in L^2 \ ,
\eeq
\beq
\label{3.62e}
x A_+ \in H^{k-1} \qquad , \qquad x \dot{A}_+ \in L^{3/2} \ ,
\eeq
$$ x A_+ \in \dot{H}^{-1/2-\mu} \ , \quad A_+\ ,  x \dot{A}_+ \in \dot{H}^{-3/2 - \mu} \ , \quad
\dot{A}_+ \in \dot{H}^{-5/2 - \mu} \ . \eqno(3.63)_{\mu}$$

\noi Let $B_{0L}$ and $h$ be defined by (1.17) and (\ref{3.31e}). \par

Then the following estimates hold~:
$$\parallel \omega^m B_{0L} \parallel_2 \ \leq C\ t^{m-1/2 + (\beta_0 - 1)(m+3/2 + \mu )} \left (
\parallel A_+ ; \dot{H}^{-3/2-\mu} \parallel \  + \ \parallel \dot{A}_+ ; \dot{H}^{-5/2 - \mu} \parallel
\right ) \eqno(3.64)$$

\noi for all $m \geq 0$,
$$\parallel \omega^m h \parallel_2 \ \leq 2t^{m-3/2} \left ( 1 \vee t^{(\beta_0 - 1)(m-1/2 +
\mu)} \right ) \left ( \parallel A_+; \dot{H}^{\rho} \parallel \ + \ \parallel\dot{A}_+; \dot{H}^{\rho -
1} \parallel \right )\eqno(3.65)$$

\noi for all $m \leq k+1$, where
$$\rho = (m-2) \vee (-3/2 - \mu) = - 3/2 - \mu + (m - 1/2 + \mu ) \vee 0 \ , \eqno(3.66)$$
$$\parallel \omega^m \partial_t h \parallel_2\ \leq C \ t^{m-5/2}  \left ( 1 \vee t^{(\beta_0 - 1)(m-1/2 +
\mu)} \right ) \left ( \parallel x A_+; \dot{H}^{\rho + 1} \parallel \right . $$
$$\left . + \ \parallel x \dot{A}_+; \dot{H}^{\rho} \parallel \ + \ \parallel A_+; \dot{H}^{\rho}
\parallel \ + \ \parallel\dot{A}_+; \dot{H}^{\rho -
1} \parallel \right )\eqno(3.67)$$

\noi for all $m \leq k$ and $\rho$ given by (3.66), 
$$\parallel h \parallel_{\infty} \ \leq C(A_+, \dot{A}_+) \eqno(3.68)$$

\noi where the constant depends on $(A_+, \dot{A}_+)$ through the norms in (\ref{3.61e})
(3.63)$_{\mu}$.} \\

\noi {\bf Proof.} (3.64) follows immediately from the definitions (1.4) (1.13) and (1.17) of $A_0$ and
$B_{0L}$, from (\ref{2.6e}) and from (3.63)$_{\mu}$. \par

In order to derive the estimates of $h$, it is convenient to come back to the variable $A_0$. The
definition (\ref{3.31e}) of $h$ can be rewritten as
$$h = 2 t^2 \omega^{-2} D_0^{-1} A_{0S} = D_0^{-1} f \eqno(3.69)$$

\noi where
$$ f = 2 \omega^{-2} A_{0S} \ , \eqno(3.70)$$  

\noi $A_{0S}$ is defined by
$$A_{0S} = t^{-1} \ D_0^{-1}\ B_{0S} = \chi_S \ A_0 \equiv F^{-1}\left ( 1 - \chi (\xi t^{1-\beta_0})
\right ) FA_0 \eqno(3.71)$$

\noi and $\chi$ is defined before (1.16). \par

(3.65). We estimate
$$\parallel \omega^m h \parallel_2 \ = t^{m-3/2} \parallel \omega^m f \parallel _2 \ = 2 t^{m-3/2} \parallel 
\omega^{m-2} A_{0S} \parallel _2$$ 
$$\leq 2 t^{m-3/2} \left ( 1 \vee t^{(\beta_0 - 1)(m-1/2 + \mu )} \right ) \parallel  \omega^{\rho} A_0 \parallel
_2 \eqno(3.72)$$

\noi and the result follows from the assumptions (\ref{3.61e}) $(3.63)_{\mu}$.\par

(3.66). We use in addition the commutation relations
$$t \partial_t = D_0^{-1} P D_0 \quad , \quad P \omega^{-j} = \omega^{-j} (P + j) \quad , \quad [P,
e^{i\omega t}] = 0 \eqno(3.73)$$

\noi where $P$ is the dilation generator
$$P = t \partial_t + x \cdot \nabla \ .$$

\noi In particular
$$\partial_t h = t^{-1} \ t \ \partial_t \ D_0^{-1} f = t^{-1}D_0^{-1}Pf \ . \eqno(3.74)$$

\noi Using the commùutation relations (3.73), we compute 
$$\begin{array}{ll}
 (1/2) Pf &=\omega^{-2} (P+2) A_{0S} \\ &\\
&= \omega^{-2} \cos \omega t \ (P+2) \chi_S A_+ + \omega^{-3} \sin \omega t \ (P+3) \chi_S \dot{A}_+ \ .
\end{array}$$

Using the fact that $P + 3 = t \partial_t + \nabla \cdot x$ and the commutation relation 
$$[P, \chi_S] = - \beta F^{-1} \xi t^{1-\beta} \cdot \nabla \chi (\xi t^{1-\beta}) F \equiv
\widetilde{\chi}$$

\noi we obtain
$$\begin{array}{ll} (1/2)Pf &= \omega^{-2} \cos \omega t \ \nabla \cdot \chi_S x A_+ + \omega^{-3} \sin
\omega t \ \nabla \cdot \chi_S x \dot{A}_+ \\ & \\ &+ \omega^{-2} \cos \omega t
\ (\widetilde{\chi} - \chi_S) A_+ + \omega^{-3} \sin \omega t \ \widetilde{\chi} \dot{A}_+ \ .
\end{array}\eqno(3.75)$$

\noi We then estimate
$$\parallel \omega^m \partial_t h \parallel _2 = t^{m-5/2} \parallel \omega^m Pf \parallel_2
\eqno(3.76)$$

\noi and we estimate the contribution of the various terms of (3.75) exactly as in the proof of (3.65),
with $(m, A_+, \dot{A}_+)$ replaced by $(m - 1, xA_+, x \dot{A}_+)$ in the first two terms, and with
$\chi_S$ replaced by $\widetilde{\chi} - \chi_S$ or by $\widetilde{\chi}$ in the last two terms. This
yields (3.67).\par

(3.68). By Lemma 2.1,
$$\parallel  h \parallel_{\infty} \ \leq C \parallel \omega^{3/2- \varepsilon}h \parallel _2^{1/2} \
\parallel \omega^{3/2 + \varepsilon} h \parallel _2^{1/2} \eqno(3.77)$$

\noi and (3.68) follows from (3.65) with $0 < \varepsilon \leq (k - 1/2) \wedge (\mu + 1)$. \par\nobreak
\hfill $\sq$ \par

We now collect the results of Lemmas 3.2, 3.3 and 3.4 in order to exhibit a set of assumptions which
imply those of Proposition 3.1\\

\noi {\bf Proposition 3.2.} {\it Let $1 < k < 3/2 < \ell$. Let $\mu$, $\lambda_0$, $\lambda$, $\beta_0$,
$\beta$ and $k_+$ satisfy the conditions
$$- 1/4 < \mu \leq 1/2 \eqno(3.78)$$
$$\lambda > 0 \quad , \qquad (1<)\lambda + k < \lambda_0 < 7/6 + 2 \mu /3 (\leq 3/2) \eqno(3.79)$$
$$0 < \beta_0 \leq \beta < 2/3 \quad , \quad \beta (\ell + 1) < \lambda_0 \eqno(3.80)$$
$$\beta_0 (1/2 - \mu ) > \lambda_0 - 1 - \mu \eqno(3.81)$$
$$\beta_0 (\mu + 5/2 ) < 2 + \mu - \lambda_0 \eqno(3.82)$$
$$k_+ \geq k + 2 \quad , \quad \beta (k_+ + 1) \geq \lambda_0 \quad , \quad \beta (\ell + 3 - k_+) < 1 \
. \eqno(3.83)$$

\noi Let $w_+ \in H^{k_+}$ and let $(A_+, \dot{A}_+)$ satisfy (\ref{2.10e}) (\ref{3.62e})
(3.63)$_{\mu}$. Let $(W, S)$ be defined by (\ref{3.25e})-(\ref{3.29e}) (\ref{3.31e}). \par
Then $(W, S)$ satisfy the estimates (\ref{3.17e}) (\ref{3.18e}) (\ref{3.19e}) (\ref{3.20e})
(\ref{3.21e}) (with $0 < \eta < 1 - 3 \beta /2$ in (\ref{3.18e})) and all the assumptions of Proposition
3.1 are satisfied}.\\

\noi {\bf Proof.} The contribution of the terms not containing $w_2$ or $B_{0L}$ in $(W, S)$ and in the
remainders are estimated by Lemma 3.2 in the same way as in I. We concentrate on the remaining terms.
The terms containing $w_2$ are estimated by Lemma 3.3 in terms of $h$, and $h$ and $B_{0L}$ are
estimated by Lemma 3.4. \par

The condition (\ref{3.17e}) restricted to $w_2 = hw_0$ follows from the fact that it holds for $h$ by
(3.65) (3.68) and trivially for $w_0$, and that it is multiplicative. Together with
(\ref{3.45e}) for $W_1$, it implies (\ref{3.45e}) for $W$. \par

We next consider $R_{1\nu}(W,S)$. Its $L^2$ norm is estimated by (\ref{3.46e}). By Lemma 3.4, it
satisfies the estimate (\ref{3.19e}) provided
$$(1 - \beta_0 ) \left ( ( 1/2 - \mu ) \vee 0 \right ) < 3/2 - \lambda_0 \eqno(3.84)$$

\noi which reduces to (3.81) for $\mu \leq 1/2$. \par

Similarly, $R_{1\nu}(W,S)$ is estimated in $\dot{H}^k$ norm by (\ref{3.47e}) and satisfies the estimate
(\ref{3.20e}) for $\delta$ sufficiently small under the condition (3.84) because the time decay of
(\ref{3.47e}) is worse than that of (\ref{3.46e}) at worst by a factor $t^{k+2\delta /3}$ which is
better than the allowed $t^{\lambda_0 - \lambda}$ for $0 < 2 \delta /3 \leq \lambda_0 - \lambda -
k$.\par

We now turn to $R_{2\nu}(W,s)$. The contribution of $B_{0L}$ is estimated by (3.64) and satisfies the
estimate (\ref{3.21e}) provided 
$$m + 1/2 + (\beta_0 - 1) (m + 7/2 + \mu ) \leq - 1 - \lambda_0 + \beta (m + 1) \eqno(3.85)$$

\noi which is implied by (3.82) for $\beta_0 \leq \beta$. \par

The term containing $w_2$ is estimated by (\ref{3.48e}) and satisfies the estimate (\ref{3.21e}) by
(3.65) under the condition (3.84).\par

We remark here that the upper bound on $\lambda_0$ in (3.79) is the compatibility condition of (3.81)
(3.82). The remaining conditions in (3.78)-(3.83) come from I. \par \nobreak \hfill $\sq$ \par
  
We now comment briefly on the various parameters that occur in Proposition 3.2 and on the conditions
(3.78)-(3.83) that they have to satisfy. The parameters $k$ and $\ell$ characterize the regularity of the
spaces of resolution for $(w,s)$. As a consequence, $k$ also characterizes the regularity of $(A_+,
\dot{A}_+)$ as given by (\ref{2.10e}). The parameter $\mu$ characterizes the behaviour of
$(\widehat{A}_+, \widehat{\dot{A}}_+) = (FA_+, F\dot{A}_+)$ at $\xi = 0$ through the condition (3.63)$_{\mu}$. The
parameters $\lambda_0$ and $\lambda$ are the time decay exponents of the norms of $q$ in $L^2$ and in
$\dot{H}^k$. The $\mu$ dependent upper bound on $\lambda_0$ in (3.79) ranges over $(1, 3/2]$ when $\mu$
ranges over $(- 1/4, 1/2]$. Since the condition (3.84) saturates at $\lambda_0 < 3/2$ for $\mu \geq
1/2$, there is no point in considering values of $\mu > 1/2$. The parameters $\beta_0$ and $\beta$
characterize the splitting of $B_0$ and $B_1$ respectively into short range and long range parts, and
therefore the splitting of the Schr\"odinger equation into transport and Hamilton-Jacobi equations.
The parameter $\beta$ should not be too large and can be taken equal to $\beta_0$. The parameter
$\beta_0$ satisfies two inequalities (3.81) and (3.82) in opposite directions, depending on $\lambda_0$
and $\mu$, and expressing the fact that $B_{0S}$ and $B_{0L}$ are not too large. The upper bound on
$\lambda_0$ in (3.79) is the compatibility condition of those inequalities. Whenever it is satisfied,
the value $\beta_0 = 1/3$ is allowed. Actually both (3.81) and (3.82) reduce to that upper bound for $\beta_0 =
1/3$. Finally $k_+$ characterizes the regularity of $w_+$ and should be sufficiently large, depending
on $k$, $\ell$, $\lambda_0$ and $\beta$. \\

\noi {\bf Remark 3.1.} For $\mu = 1/2$ the short range restriction is no longer needed in the estimates
of $h$ and $\partial_t h$ in Lemma 3.4, and therefore the splitting of $B_0$ into short range and long
range parts is no longer needed, namely $B_0$ can be kept entirely in the equation for $w$.\par

We finally discuss the condition (3.63)$_{\mu}$ of Lemma 3.4. That condition restricts the behaviour
of the relevant functions for small $|\xi |$ in Fourier transformed variables. Let $A$ be any of the
functions $A_+$, $\dot{A}_+$, $xA_+$, $x\dot{A}_+$ and define $A_<$ by $\widehat{A}_< (\xi ) = \chi (\xi
) \widehat{A}(\xi )$. Then the conditions on $A_<$ contained in (3.63)$_{\mu}$ all take the form
$$A_< \in \dot{H}^{-3/2 - \nu} \eqno(3.86)$$

\noi for $\nu = \mu$, $\mu \pm 1$. We first remark that in the proof of Lemma 3.4, all such conditions
can be replaced by
$$|\xi |^{-\nu} \widehat{A}_< \in L^{\infty} \eqno(3.87)$$

\noi at the expense of inserting an additional factor $(\ell n \ t)^{1/2}$ in (3.65) in the case of
equality, namely for $m = 1/2 - \mu$. This follows from the fact that
$$\parallel |\xi |^{m-3/2} \ \widehat{A}_{S<} \parallel_2 \ \leq \left \{ \begin{array}{ll} C
\parallel | \xi |^{-\nu} \ \widehat{A}_<\parallel_{\infty} \left (1 \vee t^{(\beta_0 - 1) (m + \nu
)} \right ) &\hbox{for} \ m \not= - \nu \\ & \\  C \parallel | \xi |^{-\nu} \
\widehat{A}_<\parallel_{\infty} (\ell n \ t )^{1/2} &\hbox{for} \ m = - \nu \ .\end{array} \right .
\eqno(3.88)$$

\noi The occurrence of the factor $(\ell n \ t)^{1/2}$ is harmless for the applications. The
condition (3.86) is weaker than (3.87) as regards the behaviour of $\widehat{A}_<$ away from zero,
since it requires only $\widehat{A}_< \in L_{loc}^2$ instead of $\widehat{A}_< \in L_{loc}^{\infty}$.
Furthermore (3.86) almost follows from (3.87), up to a change of $\nu$ into $\nu + \varepsilon$. In
fact 
$$\parallel |\xi |^{-3/2-\nu} \ \widehat{A}_< \parallel_2 \ \leq C\ \varepsilon^{-1/2} \ \parallel |\xi
|^{-\nu - \varepsilon} \ \widehat{A}_<\parallel_{\infty} \eqno(3.89)$$

\noi for $\varepsilon > 0$. In addition, under the short range condition   
$$\parallel |\xi |^{-3/2-\nu} \ \widehat{A}_{S<} \parallel_2 \ \leq C (\ell n \ t )^{1/2} \ \parallel
|\xi |^{-\nu} \ \widehat{A}_<\parallel_{\infty} \eqno(3.90)$$

\noi which is the special case $m = - \nu$ of (3.88).\par

The restrictions on $(\widehat{A}_+,\widehat{\dot{A}}_+)$ at $\xi = 0$ expressed by (3.86) have the
unpleasant feature that for $\nu \geq 0$ they cannot be ensured by imposing decay of $(A_+, \dot{A}_+)$ at
infinity in space and that they require in addition some moment conditions. For instance even for $A
\in {\cal S}$ one has
$$\omega^{-3/2-\nu} \ A_< = C |x|^{-3/2+ \nu} \ * \ A_<$$

\noi for $|\nu| < 3/2$ \cite{9r}, which behaves as 
$$|x|^{-3/2+ \nu} \ \int A\ dx$$

\noi as $|x| \to \infty$ and therefore cannot be in $L^2$ for $\nu \geq 0$ unless $\int A\ dx = 0$.
More generally when $\nu$ increases, vanishing of the $n$-th moment of $A$ is necessary as soon as
$\nu \geq n$. Actually the parameter $\mu$ in (3.63)$_{\mu}$ has been introduced in order to minimize
the number of such conditions by taking $\mu$ small. \par

We now give sufficient conditions on $(A_+, \dot{A}_+)$ in terms of space decay and vanishing of
suitable moments so as to ensure the low frequency part of $(3.63)_{\mu}$. \\

\noi {\bf Lemma 3.5.} {\it Let $-1 < \mu < 1$. Let $(A_+, \dot{A}_+)$ satisfy (\ref{3.61e}) (\ref{3.62e})
and in addition
$$x\ A_+ \in L^{3/(2 + \mu) \vee 2}\ , \ \int \dot{A}_+ \ dx = 0 \ , \ <x>^{1 + \mu + \varepsilon} \
\dot{A}_+ \in L^1 \ , \eqno(3.91)$$
$$A_+, x \dot{A}_+ \in L^{3/(3 + \mu )} \qquad \hbox{for} \ \mu < 0 \ , \eqno(3.92)$$
$$\int A_+ \ dx = \int x \ \dot{A}_+ \ dx = 0 \quad , \quad <x>^{\mu + \varepsilon} \ A_+ \in L^1 \qquad
\hbox{for} \ \mu \geq 0 \ . \eqno(3.93)$$

\noi Then $(3.63)_{\mu}$ holds.} \\

\noi {\bf Proof.} The high frequency part of $(A_+, \dot{A}_+)$ is controlled by (\ref{3.61e})
(\ref{3.62e}) and it is sufficient to consider $(A_{+<}, \dot{A}_{+<})$, although in some cases the
high frequency parts are also controlled by (3.91) (3.92).\par

We first consider $xA_+$. For $- 1/2 \leq \mu < 1$, we estimate 
$$\parallel \omega^{-1/2-\mu} \ x \ A_+ \parallel_2 \ \leq C \parallel x\ A_+ \parallel_{3/(2 + \mu)}
\eqno(3.94)$$ 

\noi by Lemma 2.1. For $\mu \leq - 1/2$, we estimate simply
$$\parallel \omega^{-1/2-\mu} \ x\ A_{+<} \parallel_2 \ \leq C \parallel x\ A_+ \parallel_2 \ .
\eqno(3.95)$$

\noi We next consider $A_+$ and $x\dot{A}_+$ together and we use $A$ to denote either of them. For $- 1
< \mu < 0$, we estimate
$$\parallel \omega^{-3/2-\mu} A \parallel_2 \ \leq C \parallel A \parallel_{3/(3 + \mu )} \eqno(3.96)$$

\noi by Lemma 2.1. For $\mu \geq 0$, we estimate
$$|\xi|^{- \mu - \varepsilon} |\widehat{A}(\xi )| = (2 \pi)^{-3/2} \ |\xi |^{-\mu - \varepsilon} \left
| \int dx \left ( \exp (-ix \xi) - 1\right ) A(x) \right |$$ $$ \leq C \parallel |x|^{\mu +
\varepsilon}\ A \parallel_1 \eqno(3.97)$$

\noi for $0 \leq \mu + \varepsilon \leq 1$. The required estimate then follows from (\ref{3.61e})
(\ref{3.62e}) (3.91) (3.93).   \par

We finally consider $\dot{A}_+$. For $\mu < 0$, we apply the previous result with $A$ replaced by
$\dot{A}_+$ and $\mu$ replaced by $\mu - 1$. For $\mu \geq 0$, we estimate

$$|\xi|^{- 1 - \mu - \varepsilon} |\widehat{\dot{A}}(\xi )| = (2 \pi)^{-3/2} \ |\xi |^{-1 - \mu -
\varepsilon} \int dx \left ( \exp (-ix\xi) - 1-i x \xi \right ) \dot{A}_+(x) $$ $$ \leq C \parallel |x|^{1 + \mu
+ \varepsilon}\ \dot{A}_+ \parallel_1 \eqno(3.98)$$ 

\noi for $0 \leq \mu + \varepsilon \leq 1$. The required estimate then follows from (\ref{3.61e}) (3.91)
(3.93). \par \nobreak \hfill $\sq$ \par
  
\mysection{Wave operators and asymptotics for (u, A)} 
\hspace*{\parindent}
In this section we complete the construction of the wave operators for
the system (1.1) (1.2) and we derive asymptotic properties of solutions in their range. The
construction relies in an essential way on Propositions 3.1 and 3.2. So far we have worked with the
system (1.20) for $(w, s)$ and the first task is to reconstruct the phase $\varphi$.
Corresponding to $S = s_0 + s_1$, we define $\phi = \varphi_0 + \varphi_1$ where
\beq
\label{4.1e}
\varphi_0 = - \int_1^t dt'\ t'^{-1} \ B_L \left ( w_0(t'), w_0(t')\right ) \ ,
\eeq
\beq
\label{4.2e}
\varphi_1 = - \int_t^{\infty} dt' (2t'^2)^{-1} |s_0(t')|^2 + 2 \int_t^{\infty} dt' \ t'^{-1}
\ B_L \left ( w_0(t') , w_1(t')\right ) \ , 
\eeq

\noi so that $s_0 = \nabla \varphi_0$ and $s_1 = \nabla \varphi_1$. \par

Let now $(w, s)$ be the solution of the system (1.20) constructed in Proposition
3.1 and let $(q, \sigma ) =(w, s) - (W,S)$. We define 
\bea
\label{4.3e}
&&\psi = - \int_t^{\infty} dt' (2t'^2)^{-1} \left ( \sigma \cdot (\sigma + 2S) + s_1 \cdot (
s_1 + 2s_0)\right ) (t') \nn \\
&&+ \int_t^{\infty} dt' \ t'^{-1} \left ( B_L (q, q) + 2
B_L(W,q) + B_L (w_1, w_1) \right ) (t')
\eea

\noi which is taylored to ensure that $\nabla \psi = \sigma$, given the fact that $s_0$, $s_1$
and $\sigma$ are gradients. The integral is easily seen to converge in $\dot{H}^1$ (see I.8.4), and to
satisfy 
\beq
\label{4.4e}
\parallel \nabla \psi \parallel_2 \ = \ \parallel \sigma \parallel_2 \ \leq C\
t^{-\lambda_0} \ . 
\eeq

\noi Finally we define $\varphi = \phi + \psi$ so that $\nabla \varphi = s$, and $(w,
\varphi )$ solves the system (1.18). For more details on the reconstruction of
$\varphi$ from $s$, we refer to Section 8 of I. \par

We can now define the wave operators for the system (1.1) (1.2) as follows. We start from the
asymptotic state $(u_+, A_+, \dot{A}_+)$ for $(u, A)$. We define $w_+ = Fu_+$, we define
$B_0$ by (1.4) (1.13), namely
$$A_0 = \dot{K}(t) \ A_+ + K(t) \ \dot{A}_+ = t^{-1} \ D_0\ B_0 \ ,$$

\noi and we define $(W,S)$ by (\ref{3.25e})-(\ref{3.29e}) (\ref{3.31e}).\par

We next solve the system (1.20) with infinite initial time by Propositions 3.1 and 3.2 and we reconstruct $\varphi$ from
$s$ as explained above, namely $\varphi = \varphi_0 + \varphi_1 + \psi$ with $\varphi_0$, $\varphi_1$ and $\psi$ defined by
(\ref{4.1e}) (\ref{4.2e}) (\ref{4.3e}) with $(q, \sigma ) = (w, s) - (W,S)$. We finally substitute $(w, \varphi )$ thereby
obtained into (1.11) (1.3) thereby obtaining a solution $(u, A)$ of the system (1.1) (1.2). The wave operator is defined as
the map $\Omega : (u_+, A_+, \dot{A}_+) \to (u, A)$. \par

In order to state the regularity properties of $u$ that follow in a natural way from the
previous construction, we introduce appropriate function spaces. In addition to the operators
$M = M(t)$ and $D = D(t)$ defined by (1.8) (1.9), we introduce the operator
\beq
\label{4.5e}
J = J(t) = x + it \ \nabla \ ,
\eeq

\noi the generator of Galilei transformations. The operators $M$, $D$, $J$ satisfy the
commutation relation
\beq
\label{4.6e}
i\ M\ D \ \nabla = J\ M\ D \ . 
\eeq

\noi For any interval $I \subset [1, \infty )$ and any $k \geq 0$, we define the space 
\bea
\label{4.7e}
{\cal X}^k(I) &=& \Big \{ u : D^*M^* u \in {\cal C} (I, H^k) \Big \}\nn \\
&=& \Big \{ u : <J(t)>^k \ u \in {\cal C} (I, L^2) \Big \} 
\eea
 
\noi where $< \lambda > = (1 + \lambda^2)^{1/2}$ for any real number or self-adjoint operator
$\lambda$ and where the second equality follows from (\ref{4.6e}). \par

We now collect the information obtained for the solutions of the system (1.1) (1.2) and state
the main result of this paper as follows.\\

\noi {\bf Proposition 4.1.} {\it Let $1 < k < 3/2 < \ell$. Let $\mu$, $\lambda_0$, $\lambda$, $\beta_0$, $\beta$
and $k_+$ satisfy the conditions (3.78)-(3.83). \par

Let $u_+ \in F H^{k_+}$, let $w_+ = Fu_+$ and $a_+ = |w_+|_{k_+}$. Let $(A_+, \dot{A}_+)$ satisfy (\ref{2.10e})
(\ref{3.62e}) (3.63)$_{\mu}$. Let $(W, S)$ be defined by (\ref{3.25e})-(\ref{3.29e}) (\ref{3.31e}). Then\par

(1) There exists $T$, $1 \leq T < \infty$ such that the auxiliary system (1.20) has a unique solution $(w, s) \in {\cal
C}([T, \infty ), X^{k,\ell})$ satisfying
\beq
\label{4.8e}
\parallel w(t) - W(t) \parallel_2 \ \leq C\left ( a_+, A_+, \dot{A}_+\right ) t^{-\lambda_0} \ ,
\eeq
\beq
\label{4.9e}
\parallel \omega^k(w(t) - W(t)) \parallel_2 \ \leq C\left ( a_+, A_+, \dot{A}_+\right ) t^{-\lambda} \ ,
\eeq
\beq
\label{4.10e}
\parallel \omega^m(s(t) - S(t)) \parallel_2 \ \leq C\left ( a_+, A_+, \dot{A}_+\right ) t^{-\lambda_0 + \beta m} \ {\it for}
\ 0 \leq m \leq \ell + 1 \ , 
\eeq

\noi for all $t \geq T$, where the constants $C(a_+, A_+, \dot{A}_+)$ depend on $(A_+, \dot{A}_+)$ through the norms
associated with (\ref{2.10e}) (\ref{3.62e}) (3.63)$_{\mu}$. \par

(2) Let $\phi = \varphi_0 + \varphi_1$ be defined by (\ref{4.1e}) (\ref{4.2e}), let $\varphi = \phi + \psi$ with $\psi$
defined by (\ref{4.3e}) and $(q, \sigma ) = (w, s) - (W, S)$. Let
$$u = MD \exp (- i \varphi ) w \eqno(1.11)\equiv(4.11)$$

\noi and define $A$ by (1.3) (1.4) (1.5). Then $u \in {\cal X}^k([T, \infty ))$, $(A, \partial_tA) \in {\cal C}([T, \infty
), H^k \oplus H^{k-1})$, $(u, A)$ solves the system (1.1) (1.2) and $u$ behaves asymptotically in time as $MD \exp (-i
\phi )W$ in the sense that it satisfies the following estimates~:
$$\parallel u(t) - M(t) \ D(t) \exp (- i \phi (t)) W(t) \parallel_2 \ \leq C(a_+, A_+, \dot{A}_+)
t^{-\lambda_0} \ ,\eqno(4.12)$$
$$\parallel |J(t)|^k \left ( \exp (i \phi (t, x/t)) u(t) - M(t) \ D(t) \ W(t) \right )
\parallel_2 \ \leq C(a_+, A_+, \dot{A}_+)t^{-\lambda} \ , \eqno(4.13)$$
$$\parallel u(t)- M(t) \ D(t) \exp (- i \phi (t))\  W(t)  \parallel_r \ \leq C(a_+, A_+,
\dot{A}_+)t^{-\lambda_0 +(\lambda_0 - \lambda ) \delta (r)/k} \eqno(4.14)$$

\noi for $0 \leq \delta (r) = (3/2 - 3/r) \leq k$, for all $t \geq T$.\par

Define in addition
$$A_2 = A - A_0 - A_1 (|DW|^2) \ . \eqno(4.15)$$

\noi Then A behaves asymptotically in time as $A_0 + A_1 (|DW|^2)$ in the sense that $A_2$ satisfies the following
estimates~: 
$$\parallel A_2(t)\parallel_2 \ \leq C(a_+, A_+, \dot{A}_+) \ t^{-\lambda_0 + 1/2} \eqno(4.16)$$ 
$$\parallel \nabla A_2(t)\parallel_2\ \leq C(a_+, A_+, \dot{A}_+) \ t^{-2\lambda_0-1/2 + (\lambda_0 - \lambda )3/2k}
\eqno(4.17)$$ 
$$\parallel \omega^{2k-1/2} \ A_2(t)\parallel_2 \ \leq \ C(a_+, A_+, \dot{A}_+) \ t^{-2\lambda -2k+1}
\eqno(4.18)$$ 

\noi for all $t \geq T$. \par

(3) The solution $(u, A)$ also behaves asymptotically as $(MD \exp (-i \phi )W_1$, $A_0 + A_1 (|DW_1|)^2$ in the sense
that the estimates (4.12)-(4.14) and (4.16)-(4.18) also hold with $W$ replaced by $W_1$ (see (\ref{3.25e})).}\\

\noi {\bf Sketch of proof.} Part (1) is a restatement of the conclusions of Proposition 3.1 supplemented by (\ref{4.4e})
and follows from Propositions 3.1 and 3.2. \par

Part (2) follows from Part (1) and is proved in exactly the same way as Part
(2) of Proposition I.8.1. \par

Part (3) is proved in the same way as Part (2). It follows from the fact that the only estimates of $W$ and $q = w - W$
that are used in the proof of Part (2) are (\ref{3.45e}) which also holds for $W_1$ and (\ref{4.8e}) (\ref{4.9e}) which
also hold for $w_2$. In fact, the latter estimates hold for $h$ by Lemma 3.4, especially (3.65), under the assumptions of
Proposition 3.2 and follow therefrom for $w_2$ in a trivial way. \par \nobreak \hfill $\sq$ \par

\noi{\bf Remark 4.1.} It may seem surprising that the improved asymptotic form $W$ for $w$ does not give rise to better
asymptotic estimates than the simpler form $W_1$ in the norms (4.12)-(4.14) and (4.16)-(4.18). The reason is that the
additional term $w_2$ is small and gives rise to small contributions in terms of those norms. This does not prevent that
term to give a large contribution to the time derivative $\partial_tw$ in (1.20) through the derivative term $t^{-2}\Delta
w_2$. That contribution is essential to allow for the solution of the system (1.20) without assuming the support
condition. The same phenomenon appears in \cite{8r}.\par \vskip 1 truecm

\noi {\large\bf  Acknowledgements.} We are grateful to Dr. A. Shimomura for enlightening discussions.

\end{document}